\documentclass[AMA,LATO1COL]{WileyNJD-v2}

\articletype{RESEARCH ARTICLE}%

\received{}
\revised{}
\accepted{}

\usepackage[utf8]{inputenc}
\usepackage{graphicx,amssymb,amstext,amsmath}
\usepackage{amsxtra,amssymb,latexsym, amscd,amsthm}
\usepackage{indentfirst}
\usepackage[mathscr]{eucal}
\usepackage{graphicx}
\usepackage{hyperref}
\usepackage{url}
\usepackage{bm}
\usepackage{multirow}

\usepackage{tikz}
\usetikzlibrary{matrix,shapes,arrows,positioning,chains}
\tikzset{
block/.style={
    rectangle,
    draw,
    text width=17em,
    text centered,
},
decision/.style={
    rectangle,
    draw,
    text width=17em,
    text centered,
    rounded corners
},
cloud/.style={
    draw,
    ellipse,
    minimum height=2em
},
descr/.style={
    fill=white,
    inner sep=2.5pt
},
connector/.style={
    -latex,
    font=\scriptsize
},
rectangle connector/.style={
    connector,
    to path={(\tikztostart) -- ++(#1,0pt) \tikztonodes |- (\tikztotarget) },
    pos=0.5
},
rectangle connector/.default=-2cm,
straight connector/.style={
    connector,
    to path=--(\tikztotarget) \tikztonodes
}
}

\usepackage{subcaption}
\usepackage{float}
\usepackage{tensor}

\usepackage{mathrsfs}
\usepackage{geometry}
\usepackage{fancyhdr}

\usepackage{listings}             % Include the listings-package
\lstdefinestyle{myCustomMATLABStyle}{
  language=Matlab,
  numbers=left,
  stepnumber=1,
  numbersep=10pt,
  tabsize=4,
  showspaces=false,
  showstringspaces=false
}
\usepackage{bm}

\newcommand{\bx}{\ensuremath{{\bm{x}}}}

\newcommand{\bg}{\ensuremath{{\bm{g}}}}
\newcommand{\ba}{\ensuremath{{\bm{a}}}}

\newcommand{\dunit}{\,\ensuremath{\text{mm}^2/\text{s}}}
\newcommand{\bunit}{\,\ensuremath{\text{s/mm}^2}}
\newcommand{\tunit}{\ensuremath{\text{ms}}}

\newcommand{\lunit}{\ensuremath{\mu\text{m}}}
\newcommand{\qunit}{\ensuremath{\text{T/m}}}

\newcommand{\bbd}[1]{\ensuremath{{\bm{#1}}}}

\newcommand{\bn}{\ensuremath{{\bm{n}}}}
\newcommand{\bi}{\ensuremath{{\bold{i}}}}

\newcommand{\ben}{\begin{equation*}}
\newcommand{\een}{\end{equation*}}
\newcommand{\be} [1] {\begin{equation} \label{#1}}
\newcommand{\ee}{\end{equation}}

\newcommand{\R}{\ensuremath{\mathbb{R}}}
\newcommand{\C}{\ensuremath{\mathbb{C}}}
\newcommand{\ignore}[1]{}

\newcommand{\bnu}{\ensuremath{{\bm{\nu}}}}
\newcommand{\Dintr}{\mathcal{D}_0}
\newcommand{\bug}{\ensuremath{{\bm{u_g}}}}

\newcommand{\soutnewr}[2]{#2}
\newcommand{\soutnew}[2]{#2}
\begin{document}

\title{Practical computation of the diffusion MRI signal of realistic neurons based on Laplace eigenfunctions}

\author[1]{Jing-Rebecca Li*}

\author[1]{Try Nguyen Tran}

\author[2]{Van-Dang Nguyen} 

%\authormark{AUTHOR ONE \textsc{et al}}

\address[1]{\orgname{INRIA Saclay-Equipe DEFI,
CMAP, Ecole Polytechnique}, \orgaddress{\state{Palaiseau}, \country{France}}}

\address[2]{\orgdiv{Division of Computational Science and Technology}, \orgname{KTH Royal Institute of
Technology}, \orgaddress{ \country{Sweden}}}

\corres{*Corresponding author.  \email{jingrebecca.li@inria.fr}}

%\presentaddress{This is sample for present address text this is sample for present address text}

\abstract[Abstract]{The complex transverse water proton magnetization subject to diffusion-encoding magnetic
field gradient pulses in a heterogeneous medium such as brain tissue can be modeled by the
Bloch-Torrey partial differential equation.  The spatial integral of the solution 
of this equation \soutnew{}{in realistic geometry} provides a gold-standard reference model for the diffusion MRI signal
arising from different tissue micro-structures of interest.
A closed form representation of this reference diffusion MRI signal has been derived twenty years ago, 
called Matrix Formalism that makes explicit 
the link between the Laplace eigenvalues and eigenfunctions of the biological cell and its diffusion MRI signal.  
In addition, once the Laplace eigendecomposition has been computed and saved, 
the diffusion MRI signal can be calculated for arbitrary diffusion-encoding sequences and b-values 
at negligible additional cost.  

Up to now, this representation, though mathematically elegant,
has not been often used as a practical model of the diffusion MRI signal, due to the difficulties of 
calculating the Laplace eigendecomposition in complicated geometries.  
In this paper, we present a simulation framework that we have implemented inside the MATLAB-based diffusion 
MRI simulator SpinDoctor that efficiently computes the Matrix Formalism representation for realistic neurons
using the finite elements method.  
We show the Matrix Formalism representation requires around \soutnew{100}{a few hundred} eigenmodes to match the reference signal computed by solving the Bloch-Torrey equation when the cell geometry comes from realistic neurons.  
As expected, the number of required eigenmodes to match the reference signal increases with smaller 
diffusion time and higher b-values.
We also converted the eigenvalues to a length scale and illustrated the link between the length
scale and the oscillation frequency of the eigenmode in the cell geometry. 
\soutnew{}{We gave the transformation that 
links the Laplace eigenfunctions to the eigenfunctions of the Bloch-Torrey operator and 
computed the Bloch-Torrey eigenfunctions and eigenvalues.}
This work is another step in bringing advanced mathematical tools and numerical method development 
to the simulation and modeling of diffusion MRI.
}


%\jnlcitation{\cname{%
%\author{Williams K.}, 
%\author{B. Hoskins}, 
%\author{R. Lee}, 
%\author{G. Masato}, and 
%\author{T. Woollings}} (\cyear{2016}), 
%\ctitle{A regime analysis of Atlantic winter jet variability applied to evaluate HadGEM3-GC2}, \cjournal{Q.J.R. Meteorol. Soc.}, \cvol{2017;00:1--6}.}

\maketitle

\footnotetext{\textbf{Abbreviations:} PDE, partial differential equation; ODE, ordinary differential equation; FEM, finite elements method;  PGSE, pulsed-gradient spin echo; ADC, apparent diffusion coefficient; HADC, homogenized ADC; HARDI, high angular resolution diffusion MRI; STA, short time approximation; BT, Bloch-Torrey; MF, Matrix Formalism; MFGA, Matrix Formalism Gaussian approximation; }

Diffusion MRI is an imaging modality that can be used 
to probe the tissue micro-structure by encoding the incoherent motion of water molecules with magnetic field gradient pulses.  Incoherent motion during 
the diffusion-encoding time causes a signal attenuation 
from which the apparent diffusion coefficient (ADC), 
and possibly higher order diffusion terms,
can be calculated \cite{Hahn1950,Stejskal1965,Bihan1986}. 
For \soutnew{unrestricted}{free} diffusion, the root of the mean squared displacement of molecules is given by $\bar{x} = \sqrt{2\,dim \,\Dintr t}$,
where $dim$ is the spatial dimension, 
$\Dintr$ is the intrinsic diffusion coefficient, and $t$ is the diffusion time. 
In biological tissue, diffusion is usually hindered or restricted (for example, by cell membranes) 
and the mean square displacement is smaller than in the case of 
\soutnew{unrestricted}{free} diffusion.  This deviation from \soutnew{unrestricted}{free} diffusion can be used to infer information
about the tissue micro-structure.  

Using diffusion MRI to get tissue structural information in the brain has been the focus of much experimental and modeling work in recent years \cite{Palombo6671,Abe4608,ReveleyE2820, Novikov5088,8625389,8007217,8660413,8481703,Assaf2008,Alexander2010,Zhang2011,Zhang2012,Burcaw2015,Palombo2017a,Palombo2016,Ning2017}.  
In terms of modeling, the predominant approach up to now has been adding the contributions to the diffusion 
MRI signal from simple geometrical components and extracting model parameters of interest. 
Numerous biophysical models subdivide the tissue into compartments described by spheres, ellipsoids, cylinders, and the extra-cellular space
\cite{Assaf2008,Alexander2010,Zhang2011,Burcaw2015,Fieremans2011,Panagiotaki2012,Jespersen2007,Palombo2017a,DESANTIS2016,LEE2018b}.
Some model parameters of interest include axon diameter and orientation, neurite density,  dendrite structure, 
the volume fraction and size distribution of cylinder and sphere components and the effective diffusion coefficient or tensor of the extra-cellular space.
The need for a mathematically rigorous model of the diffusion MRI signal arising from realistic cellular structures was
re-iterated in recent review papers \cite{Novikov2019,Novikov2018,Fieremans2018}.  

There is a gold-standard reference model
of the diffusion MRI signal, it is the Bloch-Torrey partial differential equation (PDE) that describes 
the time evolution of the complex transverse water proton magnetization subject to diffusion-encoding magnetic
field gradient pulses.  \soutnew{}{For this model to be the gold-standard, we mean of course that 
it should be posed in realistic tissue and cell geometries.}  
This PDE can be posed in a heterogeneous medium containing different cell structures 
and the extra-cellular space.  The spatial integral of the solution 
of the PDE provides a \soutnew{gold-standard}{} reference model for the diffusion MRI signal
arising from different tissue micro-structures of interest.
Because of the high computational cost of solving the Bloch-Torrey equation in complicated cell geometries, 
this gold standard model has been used primarily as a "forward model" or "simulation framework", 
in which one changes the inputs parameters such as cell geometry, intrinsic diffusion coefficient, membrane permeability, 
and study the resulting changes to the MRI signal.  This is in contrast to "inverse models", which are used to robustly estimate 
the model parameters of interest from the MRI signal, the idea being that the "inverse models" have been formulated in such
a way so that the model parameters can be correlated to biological information in the imaging voxel.  
"Inverse models" include
the biophysical models cited above.  Nevertheless, given the recent availability 
of vastly powerful computational resources and computer memory, it is possible that simulation frameworks may become directly 
useful for parameter estimation in the future (for some recent works in this direction, see \cite{Palombo2016,Rensonnet2019}). 

Two main groups of simulation frameworks for
diffusion MRI are 1) using random walkers to mimic the diffusion process in a geometrical configuration; 
2)  solving the Bloch-Torrey partial differential equation using the finite elements method (FEM).  
The first group is referred to as Monte-Carlo simulations in the literature and previous works include \cite{Hughes1995, Yeh2013, Hall2009,Palombo2016,Balls2009}.  A GPU-based acceleration of Monte-Carlo simulation was proposed in \cite{Nguyen2018a,Waudby2011}. 
Some software packages using this approach include 
\begin{enumerate}
\item Camino Diffusion MRI Toolkit\cite{Hall2009}, developed at UCL (http://cmic.cs.ucl.ac.uk/camino/);
\item DIFSIM\cite{Balls2009}, developed at UC San Diego (http://csci.ucsd.edu/projects/simulation.html);
\item Diffusion Microscopist Simulator\cite{Yeh2013}, developed at Neurospin, CEA;
\item \soutnew{}{We mention also that the GPU-based Monte-Carlo simulation code described in \cite{Nguyen2018a} is available upon request from the authors.}
\end{enumerate}
Many works\cite{Jespersen2019,Veraart2019,Ianus2016, Drobnjak2011,Mercredi2018,Rensonnet2018,Palombo2016,Rensonnet2019}
on model formulation and validation for brain tissue diffusion MRI used Monte-Carlo simulations.
The second group of simulations, which up to now has been less often used by the diffusion MRI community, 
relies on solving the Bloch-Torrey equation in a geometrical configuration using the finite elements method (FEM)
\cite{Hagslatt2003, Loren2005, Moroney2013, Nguyen2014, Beltrachini2015}\soutnew{}{(an alternative
to the FEM is the 
finite difference method\cite{Chin2002})}.
Some of the recent works about FEM diffusion MRI simulations focused on improving its computational efficiency,
by using high-performance computing \cite{Nguyen2016a, Nguyen2018} for large-scale simulations on supercomputers, by discretization on manifolds for thin-layer and thin-tube media \cite{Nguyen2019}, 
by integrating with Cloud computing resources \soutnewr{such as Google Colaboratory notebooks working 
on a web browser or with Google Cloud Platform with MPI parallelization}{}\cite{NGUYEN2019106611}.  
Our previous works in neuron diffusion MRI simulations with FEM include the simulation of neuronal dendrites using a tree model \citep{Nguyen2015} and the demonstration that diffusion MRI signals reflect the cellular organization of cortical gray matter, these signals  being sensitive to cell size and the presence of large neurons such as the spindle (von Economo) neurons \cite{Wassermann2018, Menon662601}.  

In one recent paper\cite{lid}, we presented SpinDoctor, a MATLAB-based diffusion MRI
simulation toolbox that solves the Bloch-Torrey PDE using FEM and an adaptive 
time stepping method.  SpinDoctor provides a user-friendly interface to easily define cell configurations 
relevant to the brain white matter.  
Though the original version of SpinDoctor focused on the brain white matter, we also performed simulations 
of the diffusion MRI signal from a dendrite branch to compare the computational efficiency of SpinDoctor
with the  Monte-Carlo based simulations of Camino (http://cmic.cs.ucl.ac.uk/camino/).
It was shown that at equivalent accuracy,  SpinDoctor simulations of the extra-cellular space in the 
brain white matter is 100 times faster than Camino 
and SpinDoctor simulations of a neuronal dendrite tree is 400 times faster than Camino.
We refer the reader to \cite{lid} for detailed numerical validation of SpinDoctor simulations and timing comparisons with Camino.  

In another recent paper\cite{NeuronModulePreprint}, we presented a module of SpinDoctor called the Neuron Module that enables 
diffusion MRI simulations for a group of 36 pyramidal neurons and a group of 29 spindle neurons whose morphological descriptions were found in the neuron repository {\it NeuroMorpho.Org} \cite{ascoli9247}.  
The key to making accurate simulation possible is the use of high quality finite elements meshes for the neurons.  
For this, 
\soutnew{we used licensed software from ANSA-BETA CEA Systems 
to correct and improve}{we corrected and improved} 
the quality of the geometrical descriptions of the neurons, the details of which can be found in 
\cite{NeuronModulePreprint}.  After processing, we produced high 
quality finite elements meshes for the 65 neurons.  These finite elements meshes range from having 15163 nodes to 622553 nodes.  
Currently, the simulations in the Neuron Module enforce homogeneous Neumann boundary conditions, 
meaning the spin exchange across the cell membrane is assumed to be negligible.  
\soutnew{}{In \cite{NeuronModulePreprint} we performed an accuracy and computational timing study of the serial SpinDoctor finite elements simulation and a GPU implementation of the Monte-Carlo simulation\cite{Nguyen2018a}.  We showed that at equivalent accuracy, if only one gradient direction needs to be simulated, SpinDoctor is faster than GPU Monte-Carlo.  Because the GPU Monte-Carlo method is inherently parallel, if many gradient directions need to be simulated, there is a break-even point when GPU Monte-Carlo becomes faster than SpinDoctor.  In particular,
		      at equivalent accuracy, we showed the ratio between the GPU Monte-Carlo computational time and the SpinDoctor computational time for 1 gradient 
		      direction ranges from 31 to 72 for the whole neuron simulations.}

In this paper, we present a new module of SpinDoctor called the Matrix Formalism Module.
Taking the Bloch-Torrey equation as the gold-standard reference model, 
a closed form representation of the reference signal has been derived twenty years ago, 
that is based on the eigenvalues and eigenfunctions of the Laplace operator in the relevant cell geometry.
This representation frequently goes under the name of Matrix Formalism.
The version that uses the impulse approximation of the diffusion-encoding sequence is first found in \cite{Callaghan1997}
and the version that uses the piecewise constant approximation of the diffusion-encoding sequence is first found in \cite{Barzykin1999}.
There have been numerous works using Matrix Formalism in elementary geometries such as the line segment, the disk, and the 
sphere, as well as geometries which can be written as tensor products of these elementary geometries.  We cite \cite{Grebenkov2007,Ozarslan2009,Drobnjak2011a} and refer the reader to the literature surveys on the Matrix Formalism contained in those articles.
 
There are two advantages to the Matrix Formalism signal representation.  First, the analytical advantage is that 
this representation makes explicit 
the link between the Laplace eigenvalues and eigenfunctions of the biological cell and its diffusion MRI signal.  
This clear link may help in the formulation of reduced models of the diffusion MRI signal that is closer to the physics of the problem.  Second, the computational advantage is that once the Laplace eigendecomposition has been computed and saved, 
the diffusion MRI signal can be calculated for arbitrary diffusion-encoding sequences and b-values 
at negligible additional cost.  This will make it possible to use the Matrix Formalism as the inner loop of optimization procedures. 

Up to now, Matrix Formalism, as a closed form signal representation, though mathematically elegant,
has not been often used as a practical way of computing the diffusion MRI signal in complicated geometries
such as realistic neurons.  The calculation of the Laplace eigendecomposition in realistic neurons using Monte-Carlo based simulations would be essentially impossible due to computational time and memory limitations.  
Using the FEM, the eigenfunctions of the Laplace operator can be numerically computed in an efficient way 
and this is the approach we will describe in what follows. 

\section{Theory}
Consider a domain $\Omega$ that describes the geometry of a neuron.
In this paper, we neglect the effect of water exchange between the neuron and the extra-cellular space.
Thus, impermeable boundary conditions are imposed on $\Omega$.

\subsection{Bloch-Torrey PDE}

\label{PDEsofdiffusionMRI}

In diffusion MRI, a time-varying magnetic field gradient is
applied to the tissue to encode water diffusion.  
Denoting the effective 
time profile of the diffusion-encoding magnetic field gradient by $f(t)$, and let the
vector $\bg$ contain the amplitude and direction information of the
magnetic field gradient, the complex transverse water proton magnetization 
in the rotating frame satisfies the Bloch-Torrey PDE:
\begin{align}
\label{eq:btpde}
\begin{split}
\frac{\partial}{\partial t}{M(\bx,t)} =& -\bi\gamma f(t) \bg \cdot \bx \,M(\bx,t) + \nabla \cdot (\Dintr\nabla M(\bx,t)),\quad \bx \in \Omega,
\end{split}
\end{align}
where $\gamma=2.67513\times 10^8\,\rm rad\,s^{-1}T^{-1}$ is the
gyromagnetic ratio of the water proton, \soutnew{}{$\bi$ is the imaginary unit,}
$\Dintr$ is the intrinsic diffusion coefficient in the compartment $\Omega$. 
Neglecting water exchange between the neuron and the extra-cellular space, the 
boundary condition on the boundary \soutnew{}{$\Gamma = \partial \Omega$} is the homogeneous Neumann condition:
\begin{align}
\Dintr \nabla M(\bx,t) \cdot \bn = 0, \quad \bx \in \Gamma,
\end{align}
where $\bn$ is the unit outward pointing normal vector.
The PDE also needs initial conditions: $M(\bx,0) = \rho$,
where $\rho$ is the initial spin density.
The magnetization is a function of position $\bx$ and time $t$, 
and depends on the diffusion gradient vector $\bg$ and the time profile $f(t)$.
The Bloch-Torrey PDE is a well accepted reference model for the diffusion MRI 
signal by the diffusion MRI research community.

For simplicity, we only show results for one type of diffusion-encoding sequences in this paper:
the pulsed-gradient spin echo (PGSE) \cite{Stejskal1965} sequence.
It contains two rectangular pulses of duration $\delta$, separated by a time
interval $\Delta - \delta$, for which the effective profile $f(t)$ is
\be{eq:pgse}
f(t) =
\begin{cases}
1, \quad &t_1 \leq t \leq t_1+\delta, \\
-1,
\quad & t_1+\Delta < t \leq t_1+\Delta+\delta,\\
0, \quad & \text{otherwise,}
\end{cases}
\ee
where $t_1$ is the starting time of the first gradient pulse with $t_1
+ \Delta >T_E/2$, $T_E$ is the echo time at which the signal is measured.
The diffusion MRI signal due to spins in the domain $\Omega$ is
the space integral of $M(\bx, T_E)$ in $\Omega$:
\be{eq:signal}
S^{\text{BTPDE}} := \int_{\bx\in \Omega} M(\bx, T_E)\;d\bx.
\ee

In a diffusion MRI experiment, the pulse sequence (time profile $f(t)$) is
usually fixed, while $\bg$ is varied in amplitude (and possibly also
in direction).  The signal $S$ is usually plotted against a quantity called
the $b$-value.  The $b$-value depends on $\bg$ and $f(t)$ and is
defined as
\ben
b(\bg,f) = \gamma^2 \|\bg\|^2 \int_0^{T_E} du\left(\int_0^u f(s) ds\right)^2.
\een
The reason for these definitions is that in a homogeneous medium, the
signal attenuation is $e^{-\Dintr b}$.

An important quantity that can be derived from the diffusion MRI signal is 
the ``Apparent Diffusion Coefficient'' (ADC), defined as
(assume time profile $f$ is fixed):
\be{eq:ADCdef}
  ADC  := \left. -\frac{\partial}{\partial b} \log{\frac{S(b)}{S(0)}}\right\vert_{b=0}. 
\ee
The ADC gives an indication of the root mean squared distance travelled by water
molecules in the gradient direction
$\bug = \bg/\|\bg\|$, 
averaged over all starting positions.

\subsection{Matrix Formalism signal representation} 

It is known \soutnew{, though perhaps not as well known as it deserves to be,}that using the Matrix Formalism
\cite{Callaghan1997,Barzykin1999}, the diffusion MRI signal has the following representation for the PGSE sequence.
Let $\phi_n(\bx)$ and $\lambda_n$, $n = 1,\cdots$,  be the $L^2$-normalized eigenfunctions and eigenvalues associated to the Laplace operator 
with homogeneous Neumann boundary conditions (the surface of the neurons are supposed
impermeable):
\begin{alignat*}{3}
-\nabla \Dintr \left(\nabla \phi_n(\bx)\right) &= \lambda_n \phi_n(\bx),&\quad \bx \in \Omega,\\
\Dintr \nabla \phi_n(\bx) \cdot \bnu(\bx) &= 0, &\quad \bx \in \Gamma.
\end{alignat*}
We assume the non-negative \soutnew{}{real-valued} eigenvalues are ordered in non-decreasing order:
\ben
0 = \lambda_1 < \lambda_2 \leq \lambda_3 \leq \cdots
\een
so that $\lambda_1 = 0$ (this means the first Laplace eigenfunction is the constant function).
\soutnew{}{There is only one constant eigenfunction because} we assume the neuron is a connected domain.
Let $L$ be the diagonal matrix containing the first $N_{eig}$ Laplace eigenvalues:
\be{}
L = \text{diag} [\lambda_1,\lambda_2,\cdots, \lambda_{N_{eig}}] \in \R^{N_{eig}\times N_{eig}}.
\ee

%%%%

Let $W(\bg)$ be the ${N_{eig}\times N_{eig}}$ matrix whose entries are:
\be{}
W(\bg) =  g_x A^x+g_y A^y + g_z A^z
\ee
where \soutnew{}{$\bg = (g_x,g_y,g_z),$} and
\be{}
A^{i} =\begin{bmatrix}
a_{mn}^{i}
\end{bmatrix}, \quad i = \{x,y,z\}, \quad 1 \leq m,n \leq N_{eig},
\ee
are three symmetric ${N_{eig}\times N_{eig}}$ matrices whose entries are 
the first order moments in the coordinate directions of the product of pairs of eigenfunctions:
\ben
\begin{split}
a_{mn}^x &=  \int_{\Omega}  x\phi_m(\bx) \phi_n(\bx) d\bx,\;\\
a_{mn}^y &=  \int_{\Omega}  y\phi_m(\bx) \phi_n(\bx) d\bx,\\
a_{mn}^z &=  \int_{\Omega}  z\phi_m(\bx) \phi_n(\bx) d\bx.
\end{split}
\een

We define the complex-valued matrix $K(\bg)$: \soutnew{and diagonalize it:}{}
\be{}
K(\bg) \equiv L +\bi \gamma W(\bg).
\ee
Then the following matrix
\be{}
\label{eq:MF}
H(\bg,f) =e^{-K \delta} e^{-L (\Delta-\delta)} e^{-K^* \delta},
\ee
gives the Matrix Formalism \soutnew{signal}{}representation of the solution to the Bloch-Torrey PDE.
For a constant initial density, the Matrix Formalism signal representation 
is the entry in the first row and first column of $H(\bg,f)$, \soutnew{}{multiplied by $S_0$, 
the signal with zero b-weighting}: 
\be{}
\label{eq:sig_eig}
\soutnew{}{S^{\text{MF}}(\bg,f) = \rho |\Omega| H_{11} (\bg, f) = S_0 H_{11} (\bg, f).}
\ee
We note that in Eq. \ref{eq:MF} the matrix $L$ in the exponent is diagonal and in this case, the matrix exponential is also diagonal.  
The notation $^{*}$ denotes the matrix complex conjugate transpose. \soutnew{}{In order to calculate the non-diagonal matrix exponentials $e^{-K \delta}$ and $e^{-K^* \delta}$, we diagonalize the matrix $K(\bg)$ }
\be{}
\soutnew{}{K(\bg) =  V \Sigma V^{-1},}
\ee
\soutnew{}{where $V$ has the eigenvectors in the columns and $\Sigma$ has the eigenvalues on the diagonal.}
\soutnew{}{Then $e^{-K \delta} =V e^{-\delta \Sigma} V^{-1} $ and $e^{-K^* \delta}=\left(V^{-1}\right)^{*} e^{-\delta \Sigma^{*}}V^{*}$; hence, $H(\bg, f)$ can be written as }
\begin{equation}
\soutnew{}{ H(\bg,f) = V e^{-\delta \Sigma} V^{-1} e^{-(\Delta-\delta)L} \left(V^{-1}\right)^{*} e^{-\delta \Sigma^{*}}V^{*}. }
\end{equation}
\soutnew{}{For a derivation of the Matrix Formalism formula, see Appendix \ref{sec:appendix1}.
For the transformation between 
the Laplace eigenfunctions and the eigenfunctions of the Bloch-Torrey operator,
see Appendix \ref{sec:appendix2}}.

From the Matrix Formalism signal, the analytical expression for its ADC is the following
\soutnew{}{(for a derivation, see \cite{Haddar2018})}:
\begin{equation*}
\label{eq:DeffEigenFunc}
\frac{ADC}{\Dintr} = \sum_{n=1}^{N_{eig}} \frac{(\bug \cdot \ba_{1n})^2\lambda_n\int_0^{T_E} F(t)\left(\int_0^t e^{-\lambda_n(t-s)}f(s)ds \right)dt}{\Dintr \int_0^{T_E}F^2(t)dt},
\end{equation*} 
where \soutnew{}{$F(t) \equiv \int_0^t f(s) ds$ and } the coefficients $\ba_{1n}=\begin{bmatrix}
a_{1n}^x,  a_{1n}^y, a_{1n}^z
\end{bmatrix}^T$
are
\be{}
\begin{split}
a_{1n}^x & = \frac{1}{\sqrt{\vert\Omega\vert}}\int_{\Omega} x \phi_n(\bx) d\bx,\;\\
a_{1n}^y  &= \frac{1}{\sqrt{\vert\Omega\vert}}\int_{\Omega} y \phi_n(\bx) d\bx,\;\\
a_{1n}^z  & = \frac{1}{\sqrt{\vert\Omega\vert}}\int_{\Omega} z \phi_n(\bx) d\bx.
\end{split}
\ee
We remind the reader that the first Laplace eigenfunction is the constant function. 

To clarify the relationship between the ADC and the diffusion encoding direction $\bug$, 
we rewrite the Matrix Formalism ADC as:
\begin{equation}
\frac{ADC(\bug,f)}{\Dintr} =  \bug^T\frac{D^{\text{MF}}(f)}{\Dintr} \bug,
\end{equation} 
where the Matrix Formalism effective diffusion tensor is seen to be:
\begin{align}
\label{eq:dmf}
\frac{D^{\text{MF}}(f)}{\Dintr}= \sum_{n=1}^{N_{eig}} J(\lambda_n,f) \ba_{1n} \ba_{1n}^T, 
\end{align}
with $J(\lambda_n,f)$ depending on $\lambda_n$ and $f$:
\be{}
J(\lambda_n,f)
=\frac{\lambda_n\int_0^{T_E} F(t)\left(\int_0^t e^{-\lambda_n(t-s)}f(s)ds \right)dt}{\Dintr\int_0^{T_E}F^2(t)dt} .
\ee

We also allow the possibility of computing the Matrix Formalism Gaussian Approximation (MFGA) signal, given as
\be{}
\label{eq:sig_mfga}
S^{\text{MFGA}}(\bg,f) = \rho \vert \Omega \vert \exp\left(- \bug^TD^{\text{MF}}(f) \bug\;  b\right).
\ee

\section{Method}

In a recent work, we have taken the morphological reconstructions of some realistic neurons from NeuroMorpho.Org \cite{Ascoli2007} and processed them to create high quality finite elements meshes.  The procedure is described in \cite{NeuronModulePreprint}.
Specifically, there are 65 realistic neurons whose finite elements meshes we have made available to the public, 
with the mesh size ranges from having 15163 nodes to 622553 nodes.  These are the finite elements meshes we use
in the procedure described below to calculate the Matrix Formalism signal.
To be clear, each neuron finite elements mesh consists of
\begin{enumerate}
\item a list of $N_v$ nodes in three dimensions;
\item a list of $N_e$ tetrahedral elements ($4 \times N_e$ indices referencing the nodes);
\end{enumerate}

\subsection{Finite elements discretization of the Laplace operator}

In SpinDoctor, the finite elements space is the space of continuous piecewise linear functions on 
tetrahedral elements in three dimensions (known as $P_1$).  This space has a set of basis functions whose
number is exactly the number of finite elements nodes:
\ben
 \{\varphi_k(\bx)\}, \quad  k=1,\dots, N_v.
\een
\soutnew{}{Let the finite elements nodes be denoted $\bm{q}_1,\cdots,\bm{q}_{N_v}$.  The basis function $\varphi_k(\bx)$ is a piece-wise linear function, non-zero on the tetrahedra that 
touch the node $\bm{q}_k$, and zero on all other tetrahedra.  
On a tetrahedron that touches $\bm{q}_k$, $\varphi_k(\bx)$ is equal to $1$ on $\bm{q}_k$ and it is 
equal to $0$ on the other 3 vertices of the tetrahedron.  This completely describes the piece-wise linear
function.}

Any function in the finite elements space can be written as a linear combination of the above 
basis functions
\ben
\sum_{k=1}^{N_v} c_k \varphi_k(\bx).
\een
To discretize the Laplace operator with zero Neumann boundary conditions, two matrices are needed: $\bm M\in \R^{N_v\times N_v}$ and $\bm S\in \R^{N_v\times N_v}$, 
known in the FEM literature as the mass and stiffness matrices, respectively, which are defined as follows:
\begin{equation*}
    \bm M_{jk} =\int_\Omega \varphi_k(\bx) \varphi_j(\bx) \, d\bx, \qquad \qquad  \bm S_{jk} =\int_\Omega \Dintr\,\nabla \varphi_k(\bx) \cdot \nabla \varphi_j(\bx) \, d\bx, \quad 1 \leq j, k \leq N_v.
\end{equation*}
In SpinDoctor, these matrices are assembled from local element matrices and the assembly process is based on vectorized routines of \cite{RahmanValdman2013}, which replace expensive loops over elements by operations with 3-dimensional arrays.  All local elements matrices in the assembly of $\bm S$ and $\bm M$ are evaluated at once and stored in a full matrix of size $4 \times  4 \times N_e$, where $N_e$ denotes the number of tetrahedral elements. 

The finite elements discretization described above changes the continuous Laplace operator eigenvalue problem
to the following discrete {\it matrix} eigenvalue problem: 
\begin{equation}\label{eq:matrixeig}
\begin{split}
\text{find } \{\lambda_n, p_n\},  \; 1 \leq n \leq N_v, \; \text{ such that}\\
\lambda_n \bm M p_n =-\bm S  p_n, \; p_n\in\R^{N_v},
\end{split}
\end{equation}
where $p_n$ is the eigenvector (of $N_v$ entries) associated to the eigenvalue $\lambda_n$.
Moving back to the space of functions (the function space $P_1$), the eigenfunction $\phi_n(\bx)$ associated to the eigenvalue $\lambda_n$ is then 
\ben
\phi_n(\bx) = \sum_{j=1}^{N_v} p_n^j \varphi_j(\bx),
\een  
where the entries of the eigenvector $p_n$ are the coefficients of the finite elements basis functions.
To obtain the Matrix Formalism signal representation, we calculated 
the first moments in the three coordinate directions
of the product of pairs of eigenfunctions.  They can be written in the following form,
that involves the first moments of the finite elements basis function pairs, shown in parentheses:
\be{}
\begin{split}
a_{mn}^x &=  \sum_{j=1}^{N_v}\sum_{k=1}^{N_v}  p_n^j  p_m^k \left(\int_{\Omega} x \;\varphi_j(\bx)
 \varphi_k(\bx)d\bx\right),\\
a_{mn}^y &=  \sum_{j=1}^{N_v}\sum_{k=1}^{N_v}  p_n^j  p_m^k \left(\int_{\Omega} y \;\varphi_j(\bx)
 \varphi_k(\bx)d\bx\right),\\
a_{mn}^z &=  \sum_{j=1}^{N_v}\sum_{k=1}^{N_v}  p_n^j  p_m^k \left(\int_{\Omega} z \;\varphi_j(\bx)
 \varphi_k(\bx)d\bx\right).
\end{split}
\quad 1 \leq m,n \leq N_{eig}.
\ee

\subsection{Eigenfunction length scale}
The analytical eigenvalues of a line segment of length $H$ are 
\begin{align}
\lambda_{\{1,2,\cdots\}} &= \{0,\gamma_{l} \}, \gamma_{l} = \Dintr\left( \frac{ \pi}{H/k} \right)^2, 
 k= 1,2,\cdots
\end{align}
To make the link between the computed eigenvalue and the spatial scale of the eigenmode, 
we will convert the computed $\lambda_n$ into a length scale (from the line segment eigenvalue formula):
\begin{align}
l_s(\lambda) = \frac{ \pi}{\sqrt{\lambda/\Dintr}},
\end{align}
and characterize the computed eigenmode by ${l_s}(\lambda_n)$ instead of $\lambda_n$.
To characterize the directional contribution of the eigenmode we use the fact that 
its contribution to the ADC in the direction $\bug$
is  $ J(\lambda_n,f)(\bug \cdot \ba_{1n})^2$.  
Thus, we call $\ba_{1n}=\begin{bmatrix} a_{1n}^x,  a_{1n}^y, a_{1n}^z \end{bmatrix}^T$ the "diffusion direction" 
of the $n$th eigenmode.  
We remind that the three components of $\ba_{1n}$ are the first moments 
in the 3 principle axes directions of the associated eigenfunction.

\subsection{Eigenvalues interval and minimum length scale}

We do not want to compute the entire set of eigenvalues and eigenvectors\soutnew{}{ $\{\lambda_n, p_n\}, \ 1 \leq n\leq N_v, $}
of the matrix eigenvalue problem in 
Eq. \ref{eq:matrixeig}, because the size of $\bm M$ and $\bm S$ is 
determined by the finite elements discretization (it is equal to $N_v$, the number of finite elements nodes).  
We remind that for the 65 realistic neurons whose finite elements meshes are available in the Neuron Module
of SpinDoctor, the mesh size ranges from having 15163 nodes to 622553 nodes.  
This means most of the rapidly oscillating eigenmodes
in the matrix eigenvalue problem are linked to the finite elements discretization, and 
not the physics of the problem.  
To link with the physics of the diffusion in the cell geometry, 
we set a restricted interval in which to compute the eigenvalues.  
We set the interval to be $[0, (\pi/l_s^{min})^2 \Dintr]$, where $l_s^{min}$ is the shortest 
length scale of interest in the cell geometry.
In this way, the number of computed eigenmodes, $N_{eig}$, will be
much smaller than $N_v$.

This restricted eigenvalue interval for the matrix eigenvalue problem was implemented by called the "pdeeig" function in the MATLAB PDE Toolbox, after defining a PDE model whose PDE is the Laplace equation with
 the diffusion coefficient $\Dintr$.

\section{Results}

\label{Results}

\soutnew{All the numerical results in this Section concerns the pyramidal neuron 
{\it 02b\_pyramidal1aACC}, 
whose bounding box is $[-70,113]\lunit \times[-197,165]\lunit\times[-14,18]\lunit$.
The finite elements mesh of this neuron has 44908 nodes and 171017 tetrahedral elements.
}{}
\soutnew{}{
All the numerical results in this Section concern two neurons:}
\begin{enumerate}
\item \soutnew{}{the pyramidal neuron {\it 02b\_pyramidal1aACC},
whose bounding box is $[-70,113]\lunit \times[-197,165]\lunit\times[-14,18]\lunit$.
\soutnewr{}{The geometry is shown in Figure \ref{fig:lapeigs}.}
The finite elements mesh of this neuron has 44908 nodes and 171017 tetrahedral elements.}
\item \soutnew{}{the spindle neuron {\it 03b\_spindle4aACC}, 
whose bounding box is $[-236,102]\lunit \times[-65,68]\lunit\times[-49,6]\lunit$.
\soutnewr{}{The geometry is shown in Figure \ref{fig:BT_eig}.}
Two finite elements meshes of this neuron 
were generated.  The fine mesh has 44201 nodes and 160205 elements, the coarser mesh 
has 17370 nodes and 56163 elements.}
\end{enumerate}

The numerical solution of the Bloch-Torrey PDE was computed
using SpinDoctor to obtain the reference diffusion MRI signals.
The Bloch-Torrey PDE was discretized using P1 finite elements and solved with 
build-in MATLAB routines for ordinary differential equation (ODE) systems.  
\soutnew{The tolerances of the ODE solution of the finite elements matrix system is set to 
$atol = 10^{-4}$ (absolute tolerance) and $rtol = 10^{-2}$ (relative tolerance).}{The 
time discretization is controlled by setting the absolute and relative tolerances 
of the ODE solver.}
We refer the reader to  \cite{lid} for details on how to use SpinDoctor and 
to \cite{NeuronModulePreprint} for simulation parameters for these neuron finite elements meshes. 

\subsection{Comparing the Matrix Formalism signal with the reference signal}

The Matrix Formalism Module computes $S^{\text{MF}}$ in Eq. \ref{eq:sig_eig} for the requested b-values and diffusion-encoding sequences.  In the simulations below,  $\Dintr = 2\times 10^{-3} \dunit$.
\soutnew{}{The tolerances of the ODE solution of the finite elements matrix system is set to 
$atol = 10^{-4}$ (absolute tolerance) and $rtol = 10^{-2}$ (relative tolerance).}
We simulated two diffusion-encoding sequences:
SEQ1 (PGSE, $\delta=10.6\tunit, \Delta=13\tunit$);  
SEQ2 (PGSE, $\delta=10.6\tunit, \Delta=73\tunit$);
The set of b-values simulated are $\{0,1000,2000,3000,4000\}\bunit$.
\soutnew{}{The geometry is the pyramidal neuron {\it 02b\_pyramidal1aACC}.}

In the neuron simulations that follow, we set the minimum length scale for the 
eigenvalue problem to be $l_s^{min} = 4 \lunit$, and \soutnew{
We set} the requested \soutnew{the}{} eigenvalue interval to be $[0, (\pi/l_s^{min})^2 \Dintr]$.
In this interval, 336 Laplace eigenfunctions were found, including $\lambda_1=0$, which corresponds 
to the length scale $l_s = \infty$.
There are 6 eigenmodes with length scale $l_s > 100 \lunit$, they correspond to the length scales (rounded to the $\lunit$), 
$\{405, 343, 162, 156, 133, 127, 106\}\lunit$, respectively.

\soutnew{
In Figure \ref{fig:eig_longest} we plot the eigenfunction corresponding to the longest (finite) 
length scale, in other words, the smallest (non-zero) eigenvalue.  Its length scale is
$l_s(\lambda_n) = 405.4\lunit$ and its "diffusion direction" is parallel to 
$\ba_{1n} = [0.43, -0.9, -0.02]^T$.  This "diffusion direction" lies mostly in the $x-y$ plane
and more so in $y$ than in $x$, just like the positioning of the two main dendrite branches
of this neuron.  We conjecture that the length scale (the eigenvalue) corresponds to the "wavelength" of the 
significant oscillations of the eigenfunction in the geometry.
}{}

Next, we compare the plots of the High Angular Resolution Diffusion MRI (HARDI) signals computed in four different ways:
\begin{enumerate}
\item Reference signals from solving the Bloch-Torrey PDE \soutnew{}{($S^{\text{BTPDE}}$)}, computed in 151 diffusion directions uniformly distributed in the unit sphere;
\item  Matrix Formalism signals \soutnew{}{($S^{\text{MF}}$)} using 336 eigenfunctions found in the interval $[0, (\pi/l_s^{min})^2 \Dintr]$, $l_s^{min}=4\lunit$, computed in 151 diffusion directions uniformly distributed in the unit sphere;
\item  Matrix Formalism Gaussian Approximation signals \soutnew{}{($S^{\text{MFGA}}$)} using 336 eigenfunctions as above, computed in 151 diffusion directions uniformly distributed in the unit sphere;
\item  Matrix Formalism signals using 336 eigenfunctions as above, computed in 900 diffusion directions uniformly distributed in the unit sphere;
\end{enumerate}

\begin{figure}[ht!]
 \centering
\captionsetup[subfigure]{size=normalsize}

 \begin{subfigure}{1\textwidth}
\includegraphics[width=0.24\textwidth]{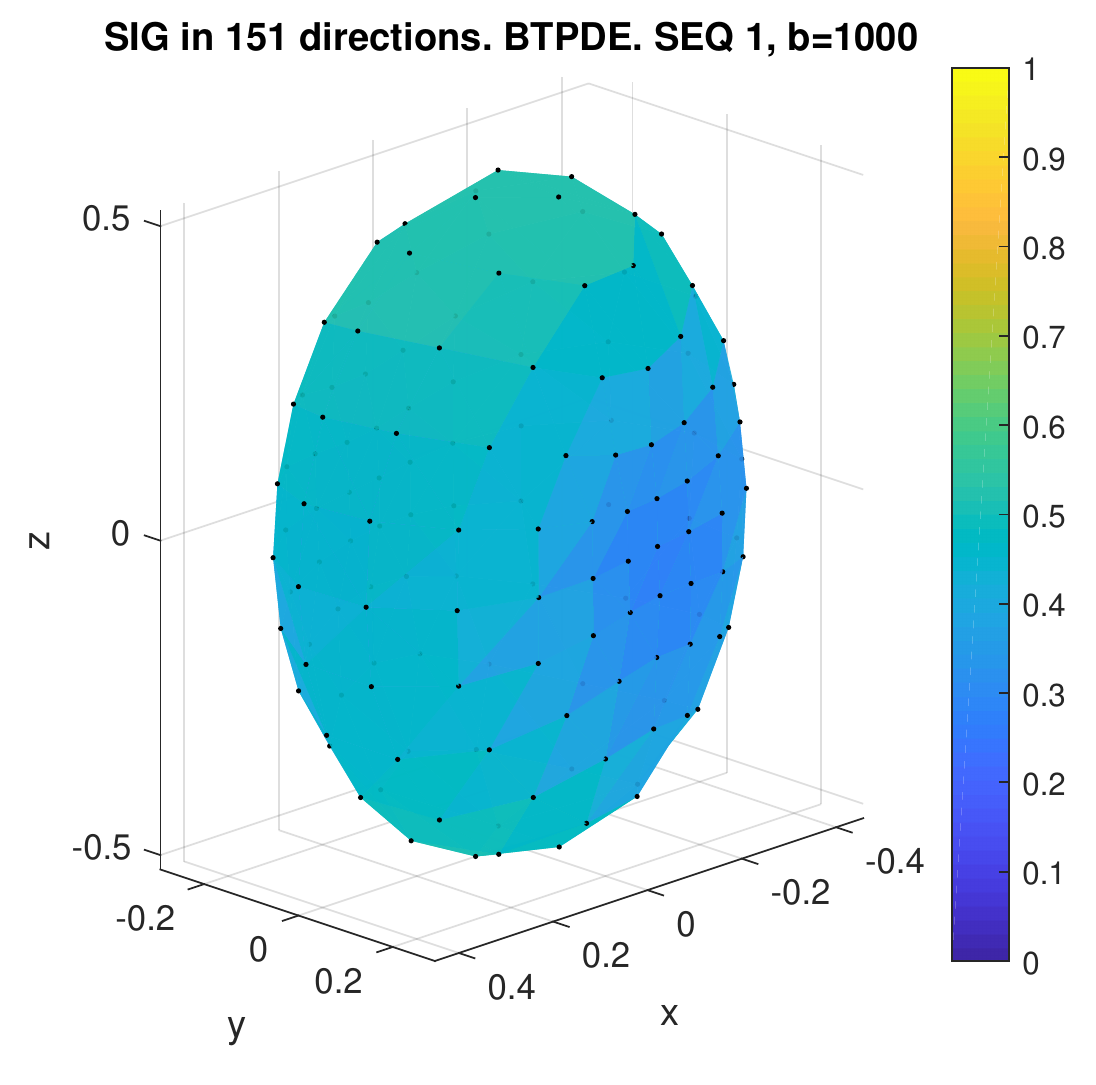}  \ 
\includegraphics[width=0.24\textwidth]{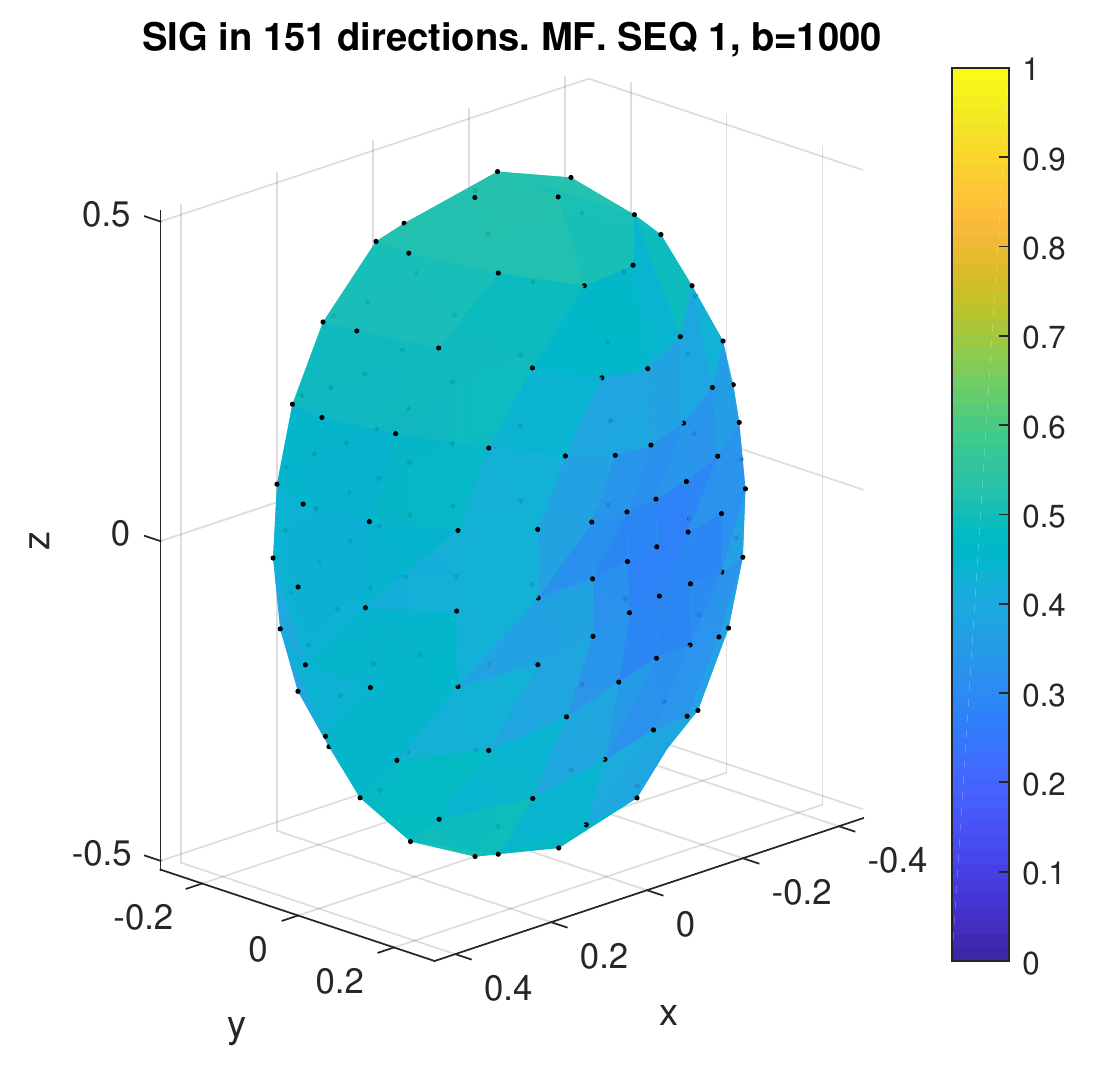} \
\includegraphics[width=0.24\textwidth]{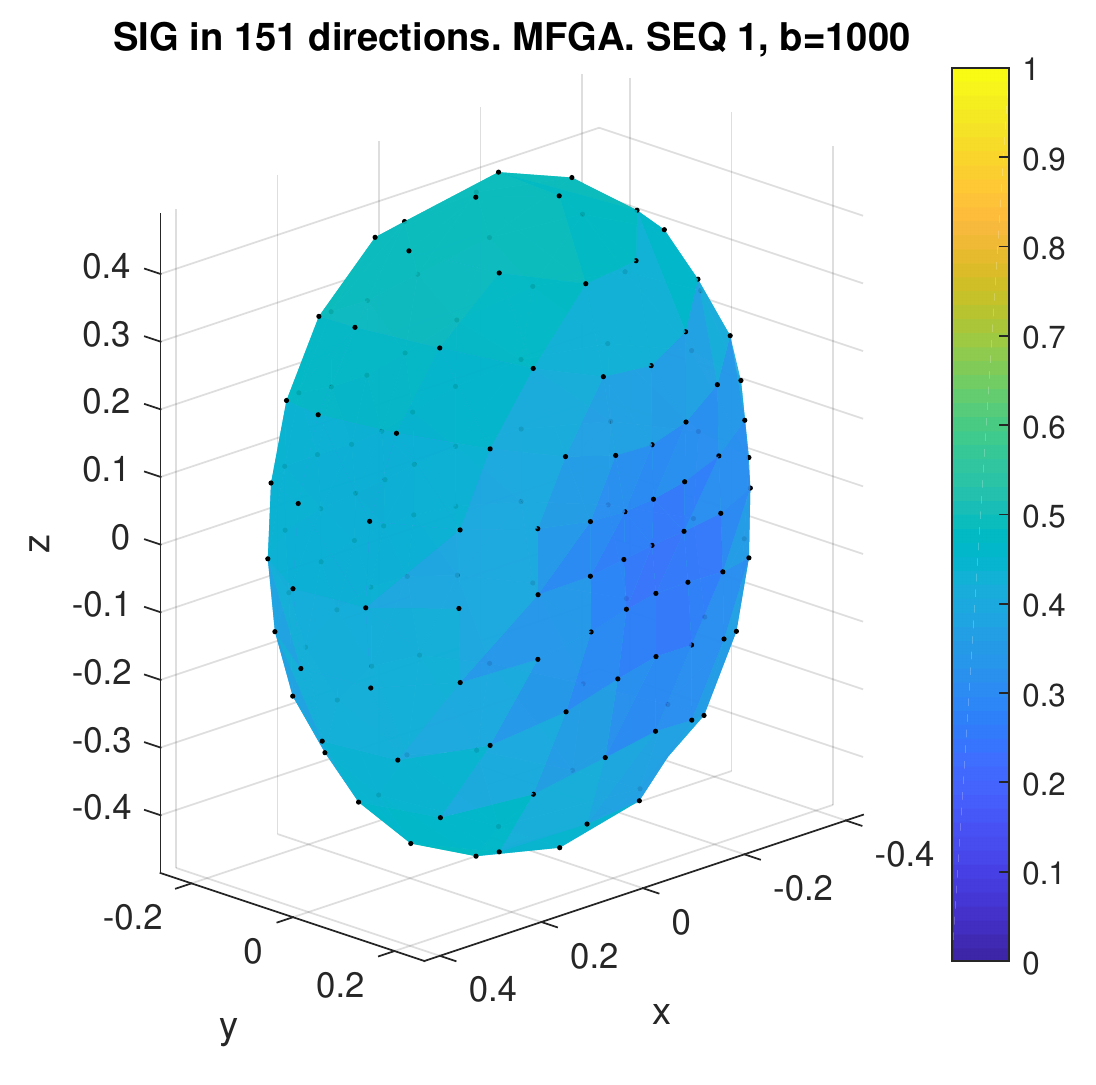}  \ 
\includegraphics[width=0.24\textwidth]{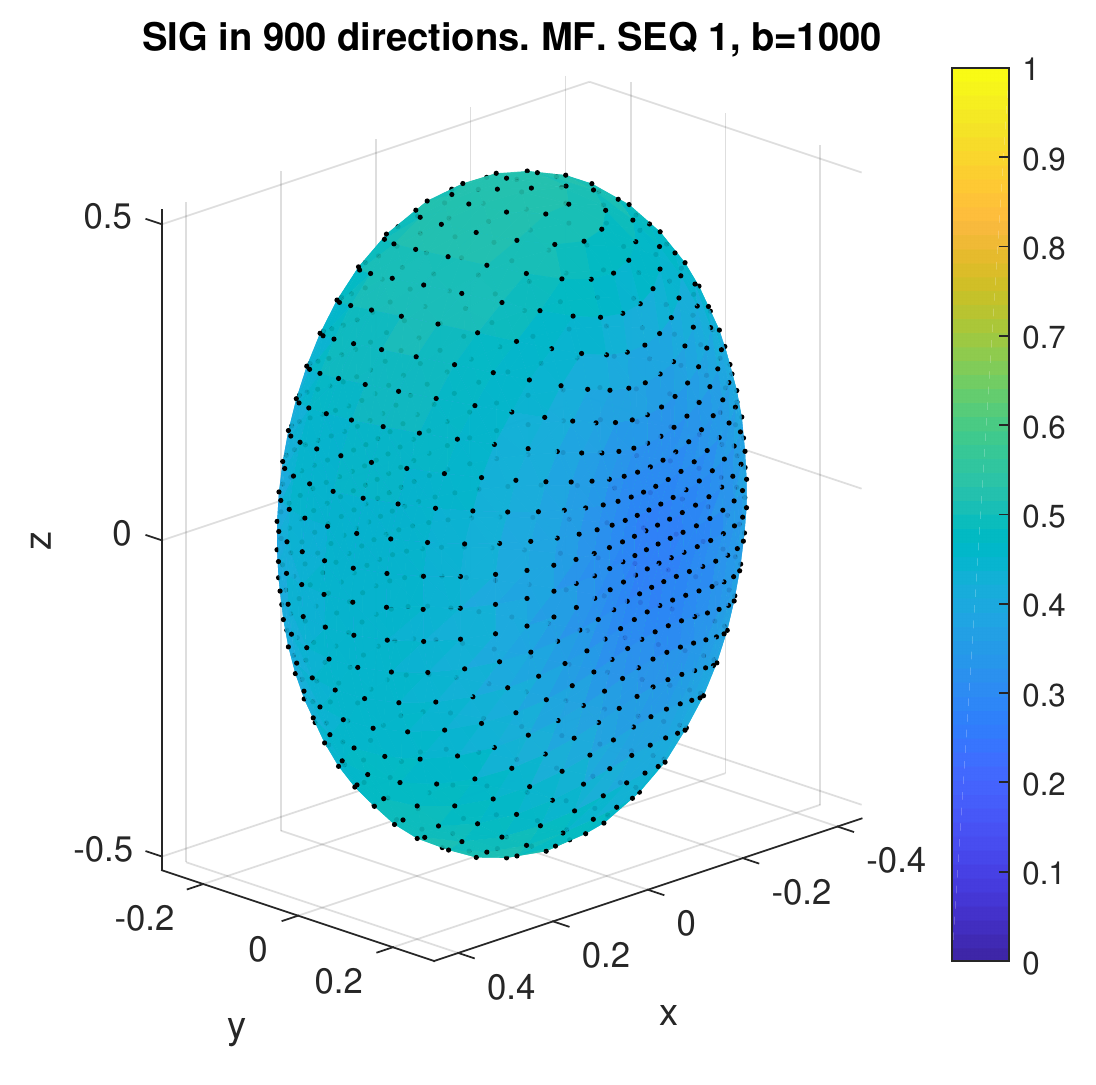}
\caption{\label{fig:sig_hardi_ex1_b1000} \soutnewr{}{SEQ1, $b = 1000\bunit$.}}
\end{subfigure}%
\vspace{0.5 cm}

 \begin{subfigure}{1\textwidth}
\includegraphics[width=0.24\textwidth]{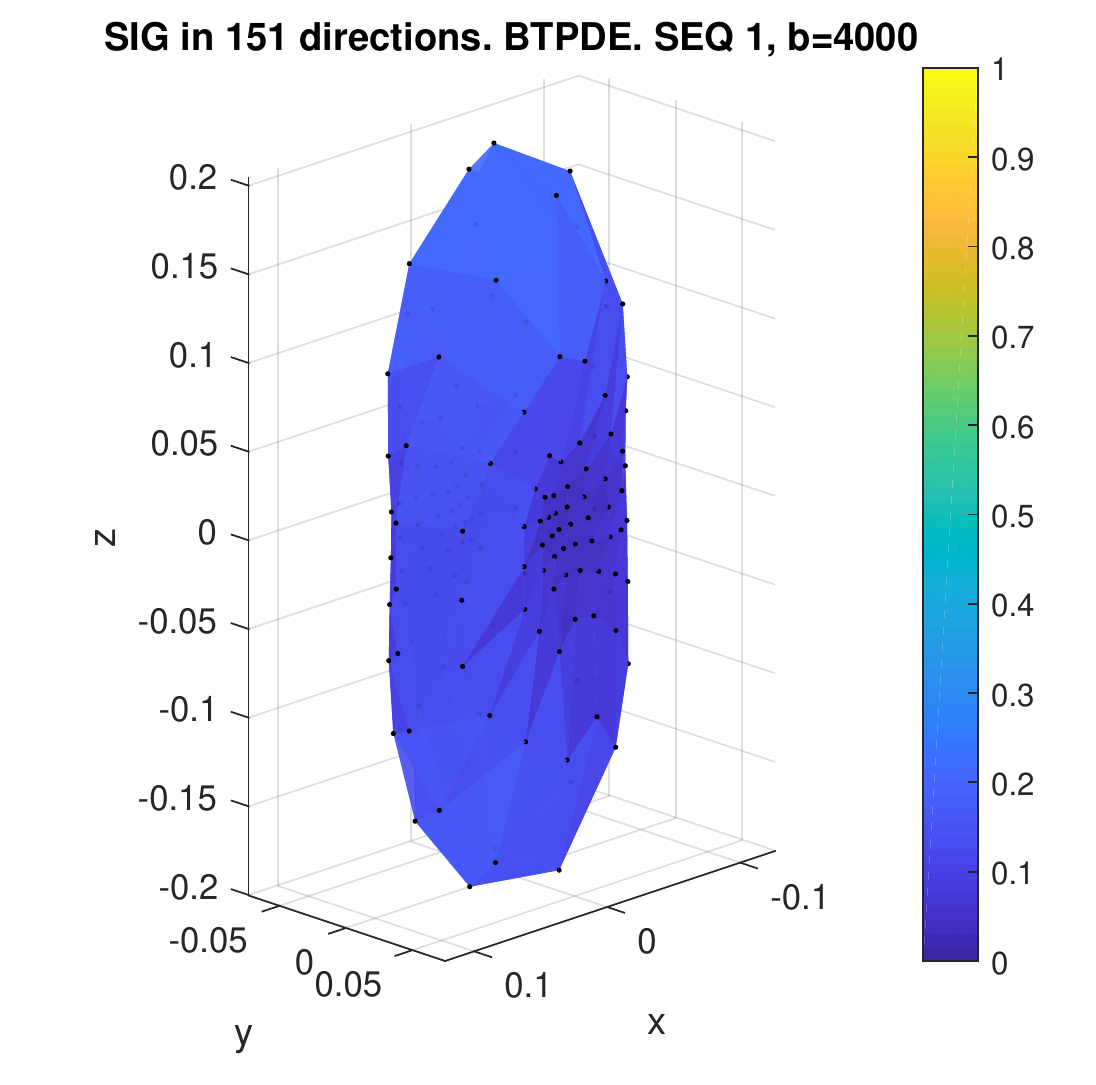}  \
\includegraphics[width=0.24\textwidth]{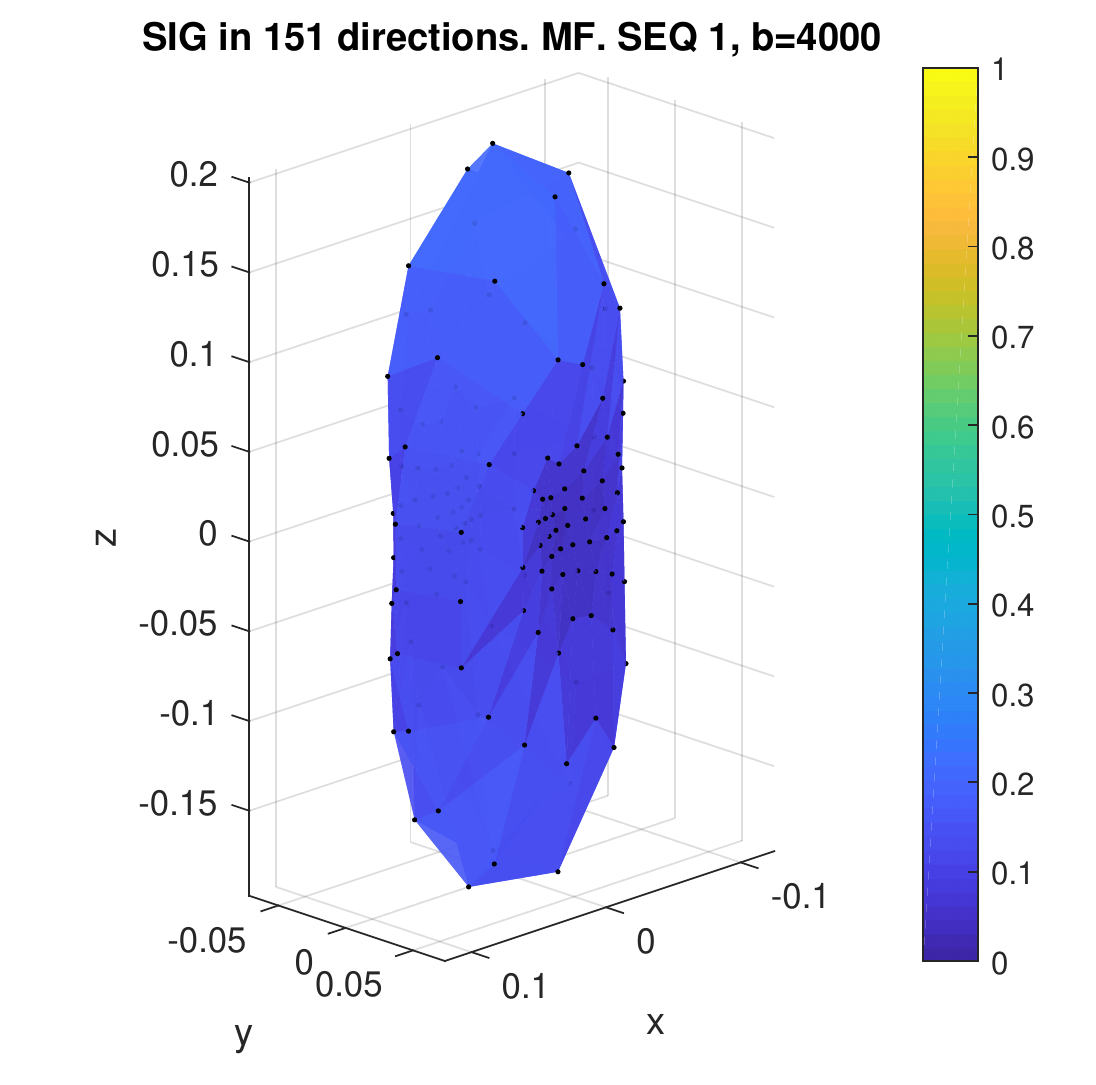} \
\includegraphics[width=0.24\textwidth]{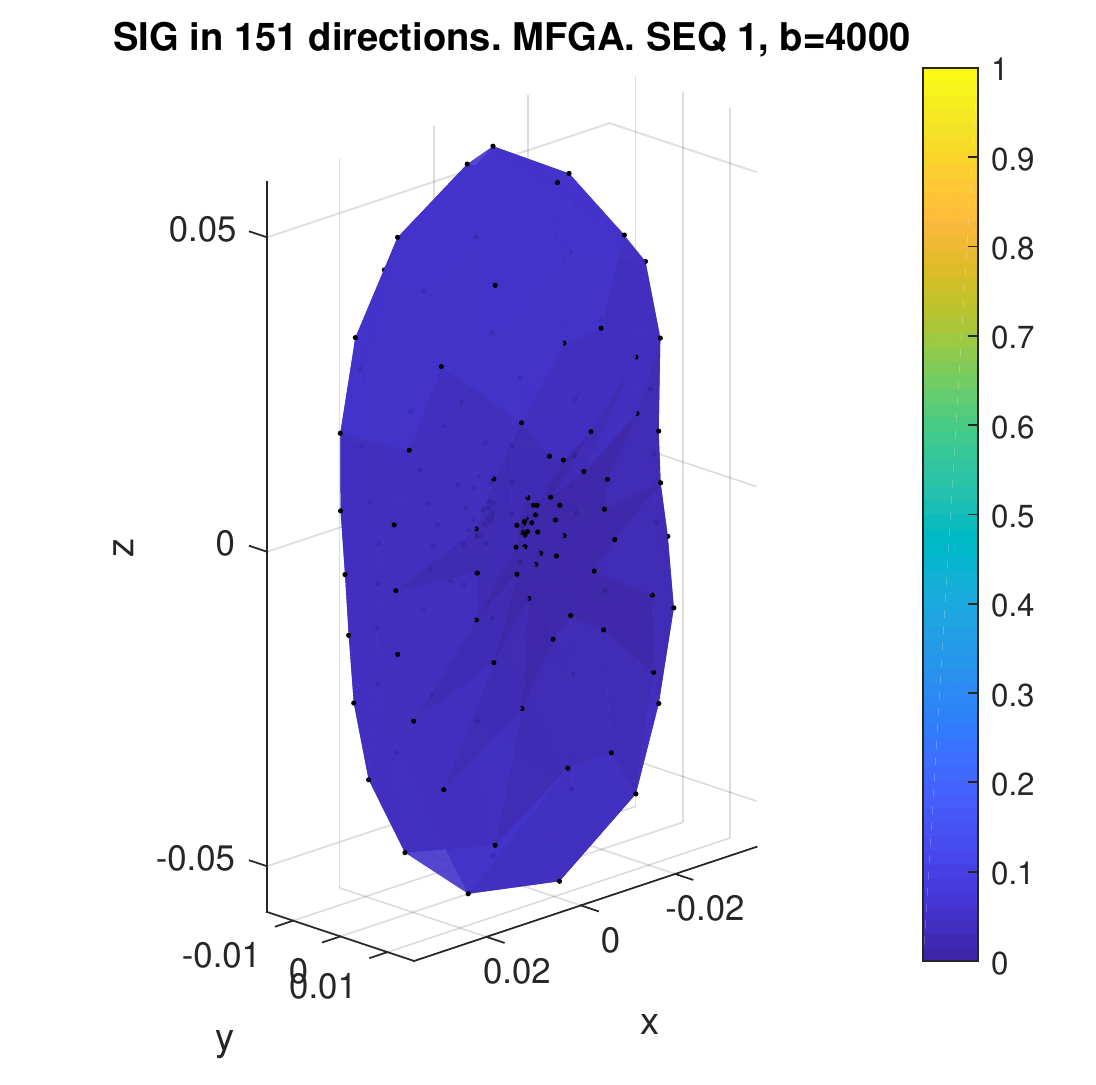}  \
\includegraphics[width=0.24\textwidth]{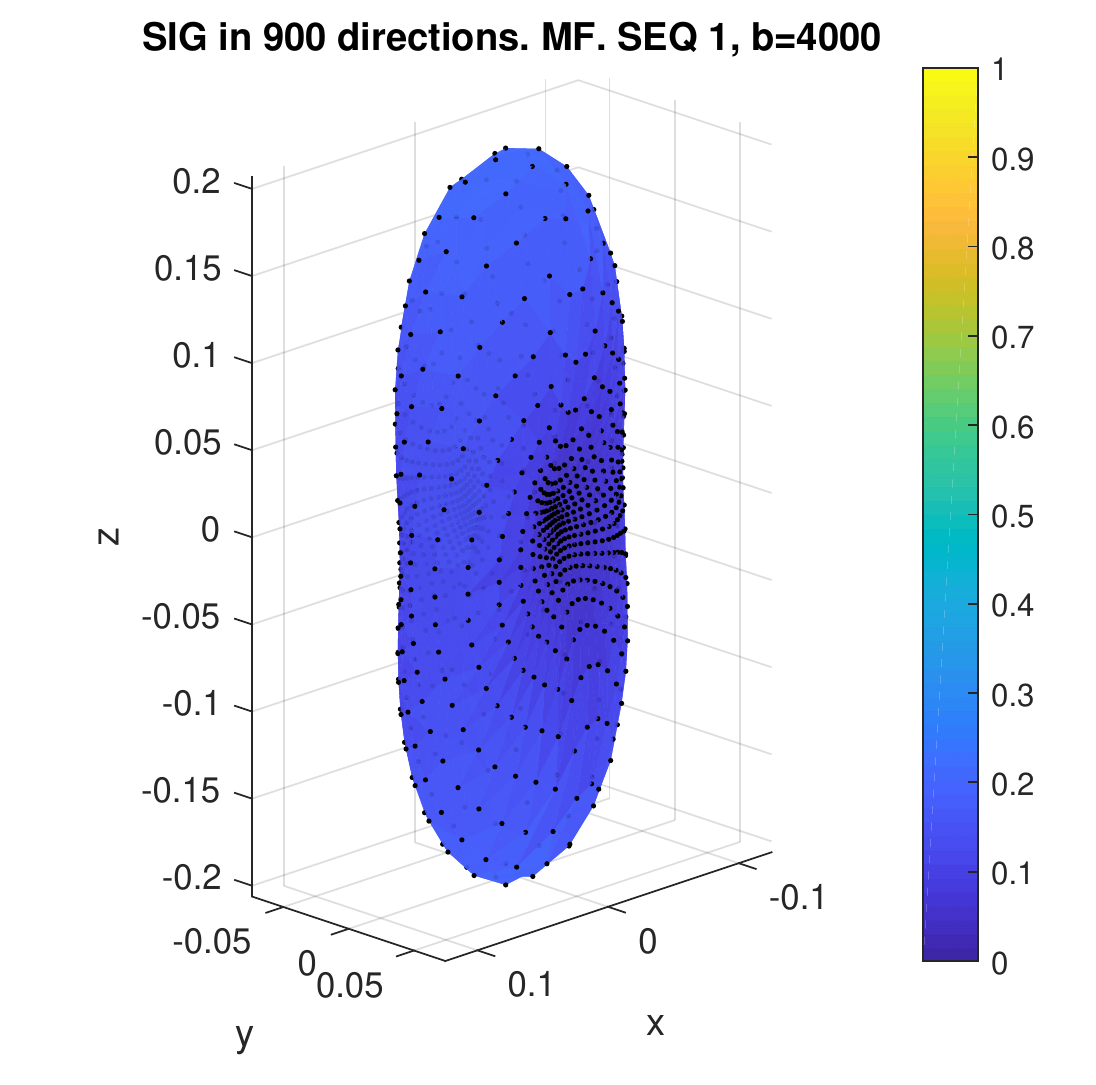}
\caption{\label{fig:sig_hardi_ex1_b4000} \soutnewr{}{SEQ1, $b = 4000\bunit$.}}
\end{subfigure}%
\vspace{0.5 cm}
 
 \begin{subfigure}{0.5\textwidth}
\includegraphics[width=0.48\textwidth]{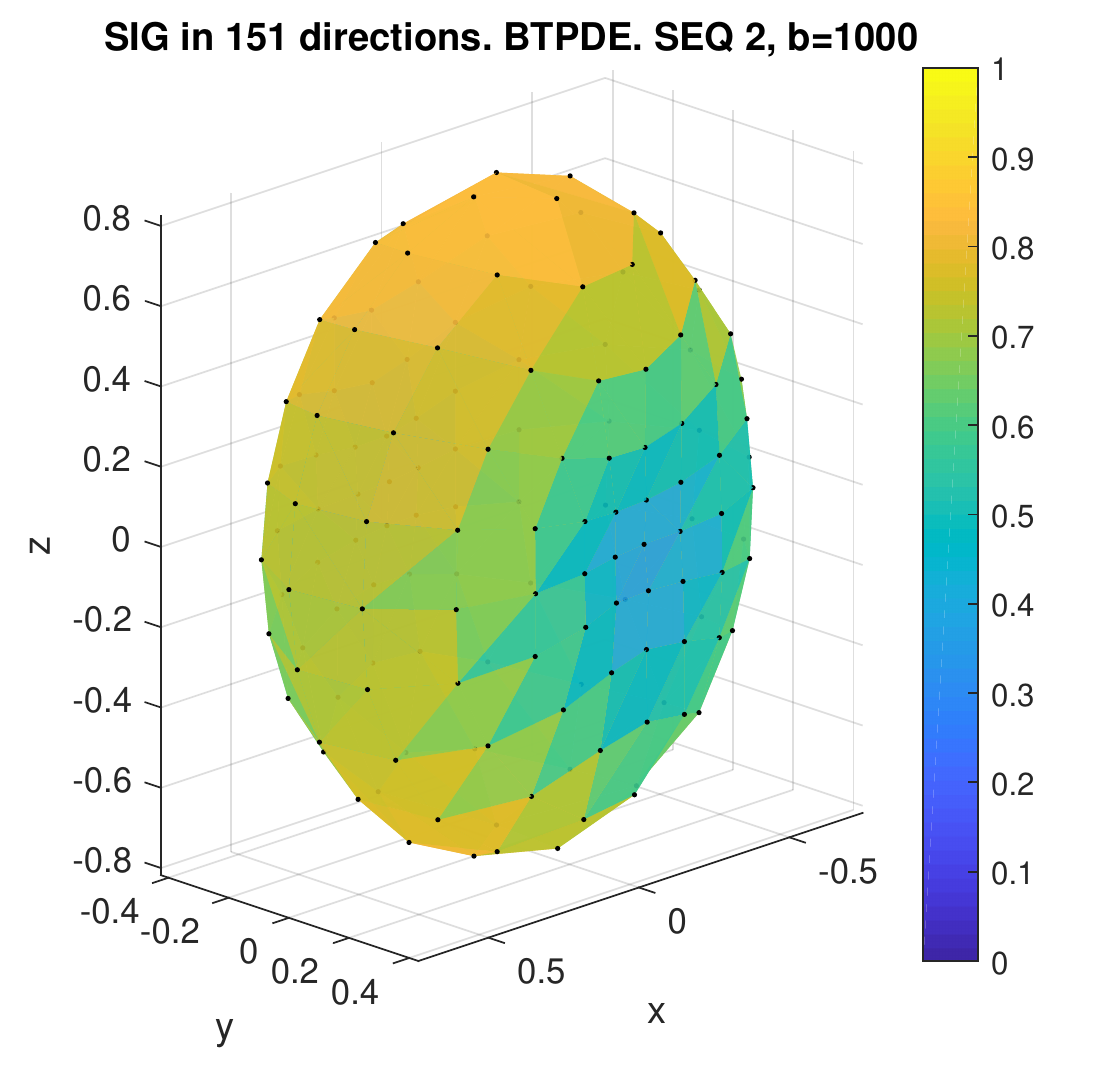}  \
\includegraphics[width=0.48\textwidth]{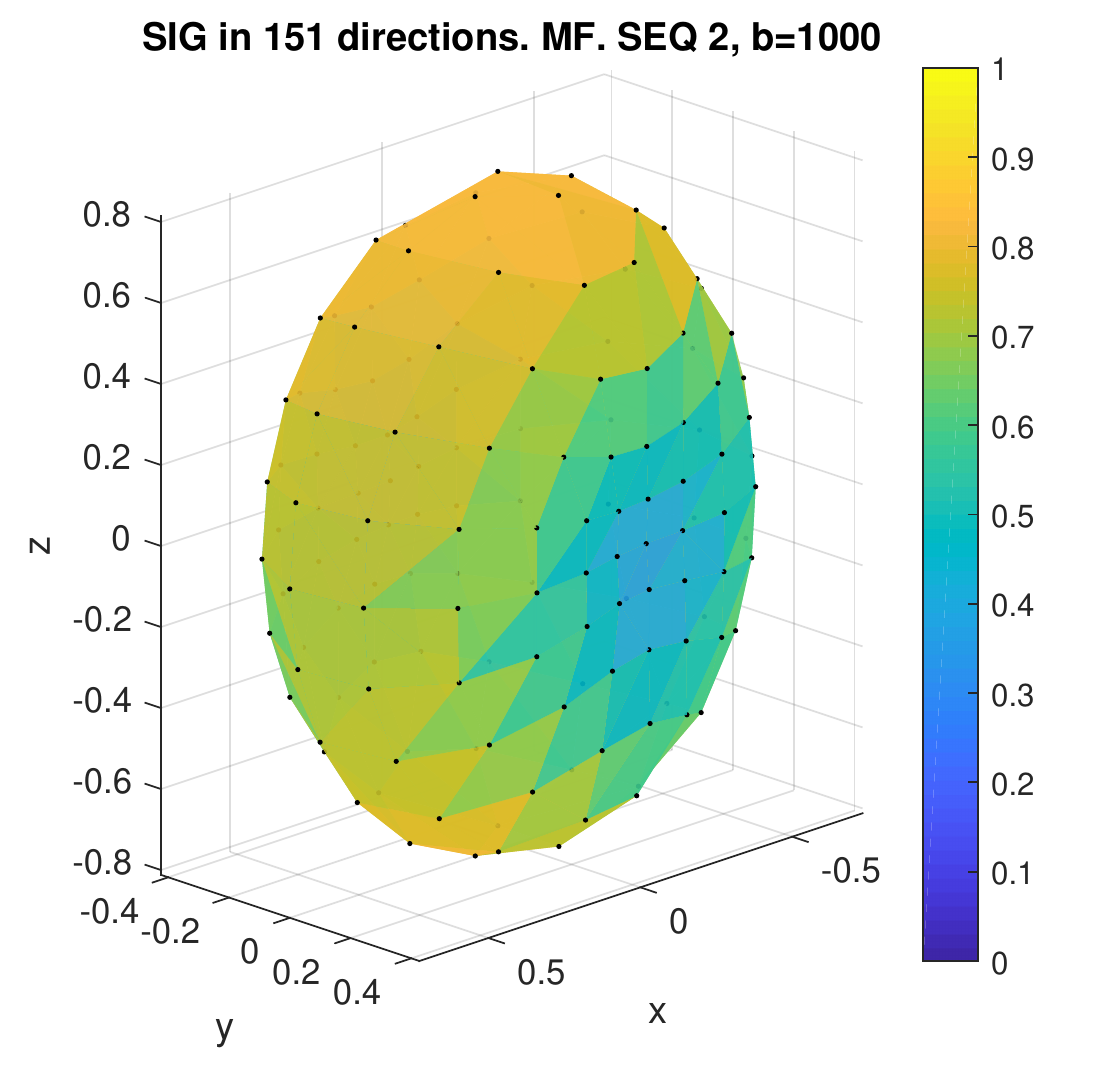} \
\caption{\label{fig:sig_hardi_ex2_b1000} \soutnewr{}{SEQ2, $b = 1000\bunit$.}}
\end{subfigure}~ \quad
 \begin{subfigure}{0.5\textwidth}
\includegraphics[width=0.48\textwidth]{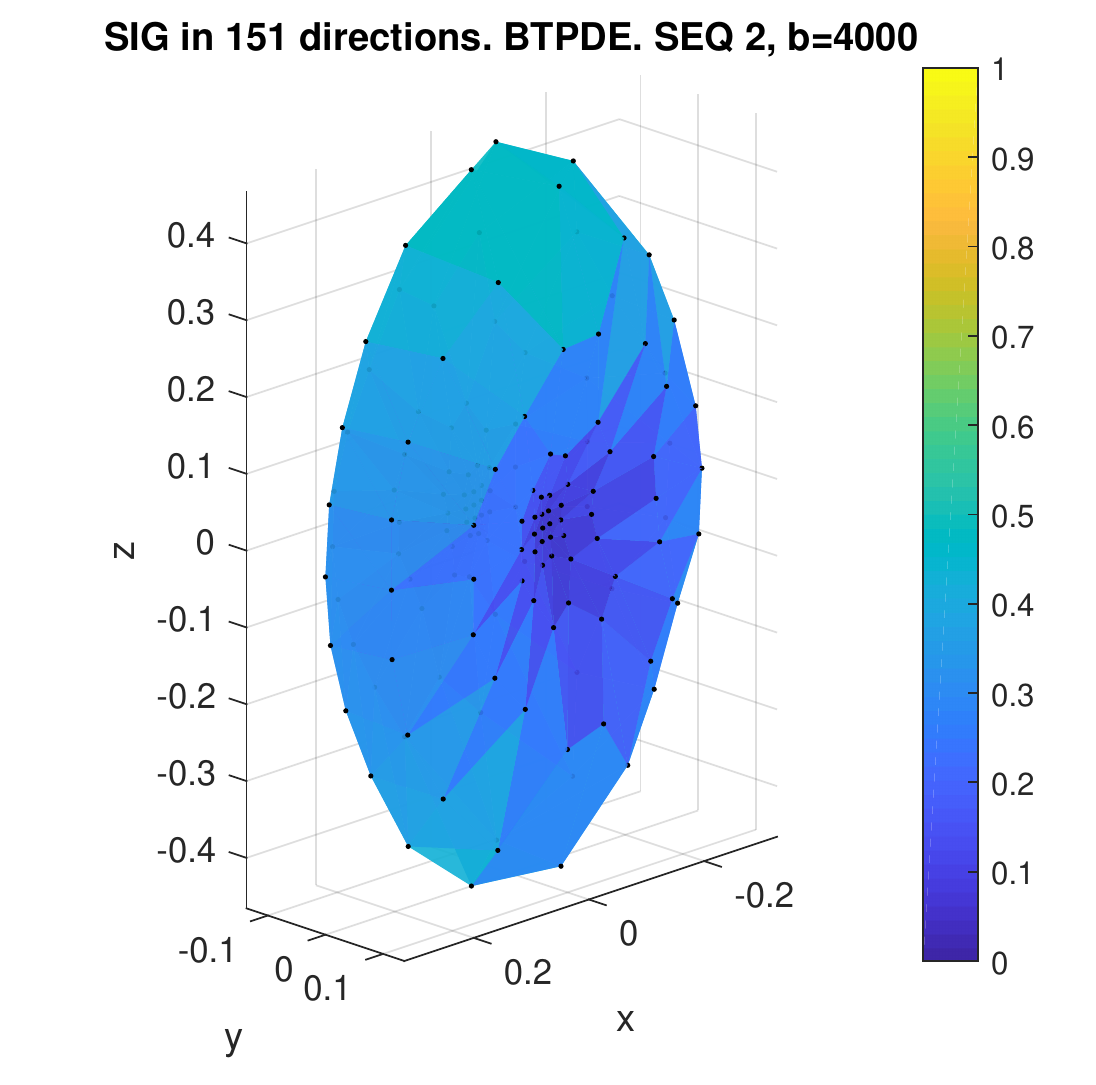}  \
\includegraphics[width=0.48\textwidth]{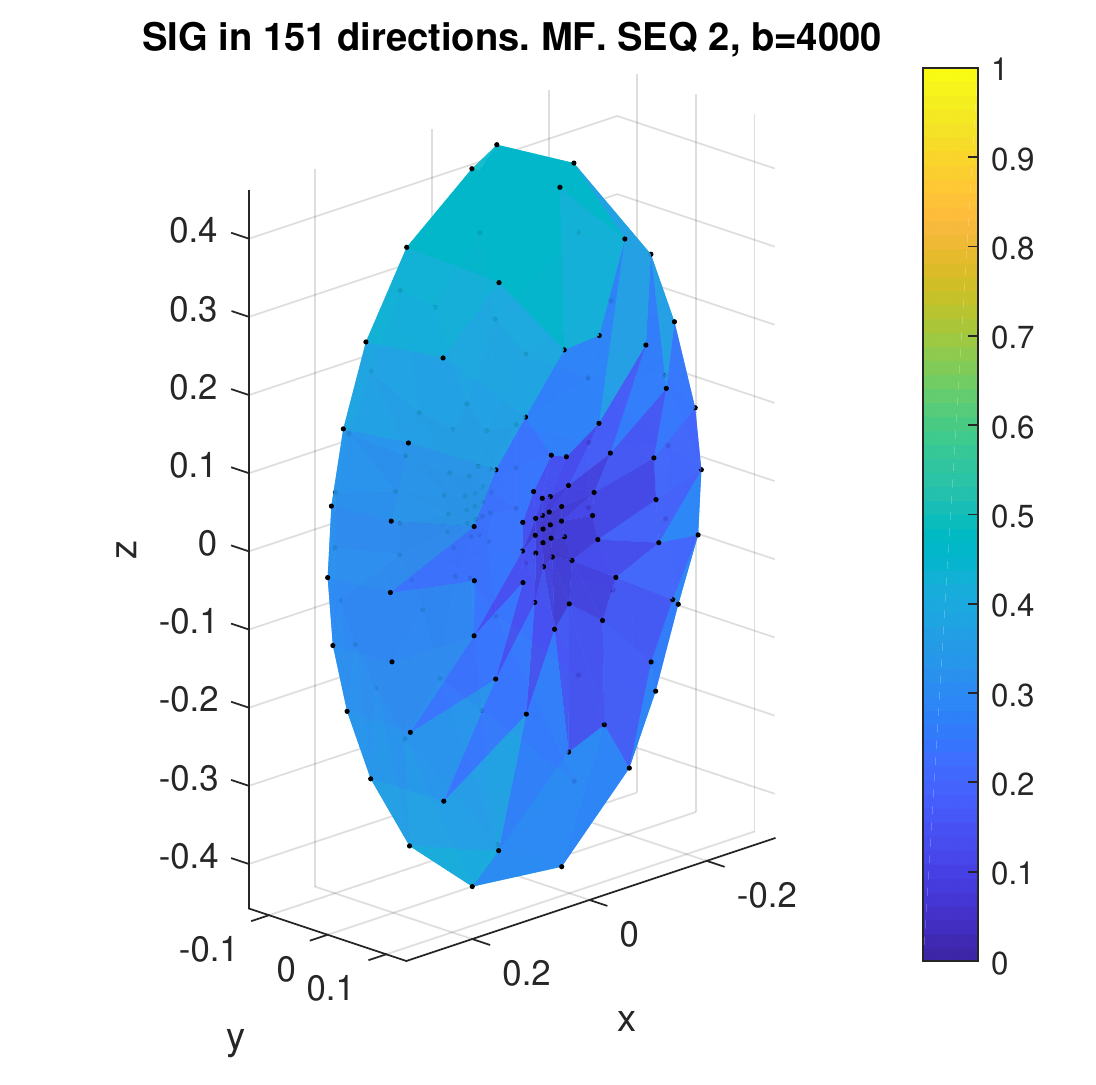} 

\caption{\label{fig:sig_hardi_ex2_b4000} \soutnewr{}{SEQ2, $b = 4000\bunit$.}}
\end{subfigure}%

\caption{ \label{fig:sigdiff} BTPDE signals, $S^{\text{BTPDE}}/S_0$, in 151 
diffusion-encoding directions.  MF signals, $S^{\text{MF}}/S_0$, in 151 diffusion-encoding directions. MFGA signals, $S^{\text{MFGA}}/S_0$, in 151 diffusion-encoding directions.  MF signals, $S^{\text{MF}}/S_0$, in 900 diffusion-encoding directions.    
The black points indicate the diffusion-encoding direction, multiplied by the magnitude of the signal attenuation.  The color indicates the value of the signal attenuation. SEQ1 is (PGSE, $\delta=10.6\tunit, \Delta=13\tunit$), 
SEQ2 is (PGSE, $\delta=10.6\tunit, \Delta=73\tunit$).
The geometry is the pyramidal neuron {\it 02b\_pyramidal1aACC}. }
\end{figure}

In Figure \ref{fig:sig_hardi_ex1_b1000} and Figure \ref{fig:sig_hardi_ex1_b4000} we show the above four HARDI signals for SEQ1. \soutnew{ at $b = 1000\bunit$ 
for SEQ1, normalized by the signal at $b=0$.  
The HARDI signal shapes are ellipsoid at this b-value and visually, 
$S^{\text{BTPDE}}$, $S^{\text{MF}}$, $S^{\text{MFGA}}$ are indistinguishable.
In Figure \ref{fig:sig_hardi_ex1_b4000} we show the four HARDI signals at $b = 4000\bunit$ 
for SEQ1, normalized by the signal at $b=0$. 
The HARDI signal shapes are no longer ellipsoid at this high b-value and visually, 
$S^{\text{BTPDE}}$, $S^{\text{MF}}$ are indistinguishable, whereas $S^{\text{MFGA}}$
is clearly different from the reference signals at this high b-value.}{At $b = 1000\bunit$ (Figure \ref{fig:sig_hardi_ex1_b1000}), 
the HARDI signal shapes are ellipsoid, and visually, 
$S^{\text{BTPDE}}$, $S^{\text{MF}}$, $S^{\text{MFGA}}$ are indistinguishable.  At 
$b = 4000\bunit$ (Figure \ref{fig:sig_hardi_ex1_b4000}),
the HARDI signal shapes are no longer ellipsoid, and visually, 
$S^{\text{BTPDE}}$, $S^{\text{MF}}$ are indistinguishable, whereas $S^{\text{MFGA}}$
is clearly different from the reference signals at this high b-value.}
\soutnew{In Figure \ref{fig:sig_hardi_ex2_b1000} we show the above four HARDI signals for SEQ2.
At $b = 1000\bunit$ 
for SEQ2, normalized by the signal at $b=0$.  
The HARDI signal shapes are ellipsoid at this b-value and visually, 
$S^{\text{BTPDE}}$, $S^{\text{MF}}$, $S^{\text{MFGA}}$ are indistinguishable.
In Figure \ref{fig:sig_hardi_ex2_b4000} we show the four HARDI signals at $b = 4000\bunit$ 
for SEQ2, normalized by the signal at $b=0$.  
The HARDI signal shapes are no longer ellipsoid at this high b-value and visually, 
$S^{\text{BTPDE}}$, $S^{\text{MF}}$ are indistinguishable, whereas $S^{\text{MFGA}}$
is clearly different from the reference signals at this high b-value.}{
For SEQ2, we see in Figure \ref{fig:sig_hardi_ex2_b1000}, that at $b = 1000\bunit$, 
the HARDI signal shapes are ellipsoid, with less signal attenuation than SEQ1. 
At $b = 4000\bunit$ (Figure \ref{fig:sig_hardi_ex2_b4000}), the HARDI signal shapes are no longer ellipsoid, however, they are 
visually quite different than the shapes for SEQ1.}

To verify numerically that the Matrix Formalism signal is close to the reference signal, 
we compute the signal differences between  $S^\text{MF}$ and the reference $S^{\text{BTPDE}}$  
over $30$ diffusion-encoding directions:
\be{}
E(f,b) = \frac{\sum_{j=1}^{30}\left( S^\text{MF}(f,\bg_j)-S^{\text{BTPDE}}(f,\bg_j)\right)^2}{\sum_{j=1}^{30}\left( S^\text{BTPDE}(f,\bg_j)\right)^2}.
\ee
The directions are uniformly distributed on the unit sphere.
The signal differences are \{1.6\%,2.2\%,0.6\%,1.9\%\} in order of 
\{(SEQ1, $b=1000\bunit$), (SEQ1, $b=4000\bunit$), (SEQ2, $b=1000\bunit$), (SEQ2, $b=4000\bunit$)\}.   
Thus, we consider the Matrix Formalism signal with the full set of 336 computed eigenvalues to be a good approximation of the 
reference signal.

\subsection{The contribution of each eigenmode to the signal}

Because $S^{\text{MF}}$ contains the contributions of all the computed eigenmodes in the requested 
interval, to get an idea of 
the importance of each eigenmode, we computed the signal difference that results when one eigenmode is removed, 
compared to using the full set of computed eigenmodes.  This signal difference is computed for each 
sequence and each b-value, averaged over 30 gradient directions.  
The directions are uniformly distributed on the unit sphere.
For the eigenfunction $i$, the signal difference is obtained as:
\be{}
E^{\text{RM}, i}(f,b) = \frac{\sum_{j=1}^{30}\left( S^\text{MF}(f,\bg_j)-S^{\text{MF}, \text{RM}, i }(f,\bg_j)\right)^2}{\sum_{j=1}^{30}\left( S^\text{MF}(f,\bg_j)\right)^2}.
\ee
The signal $S^\text{MF}$ uses the full set of computed eigenfunctions, the signal $S^{\text{MF}, \text{RM}, i }$ excludes the $i$th eigenfunction.  In the following the signal differences will be given 
for two sequences at 2 b-values, in order of 
\{(SEQ1, $b = 1000\bunit$), (SEQ1, $b=4000\bunit$), (SEQ2, $b=1000\bunit$), (SEQ2, $b=4000\bunit$)\}.  
We expect the second value to be the highest and the third value to be the lowest.
We denote the $i$th eigenfunction as "significant" if $E^{\text{RM}, i}(f,b)$ is more than a certain threshold.

\begin{figure}[ht!]
 \centering
\captionsetup[subfigure]{size=normalsize}
  \begin{subfigure}{0.48\textwidth}
\includegraphics[width=1\textwidth]{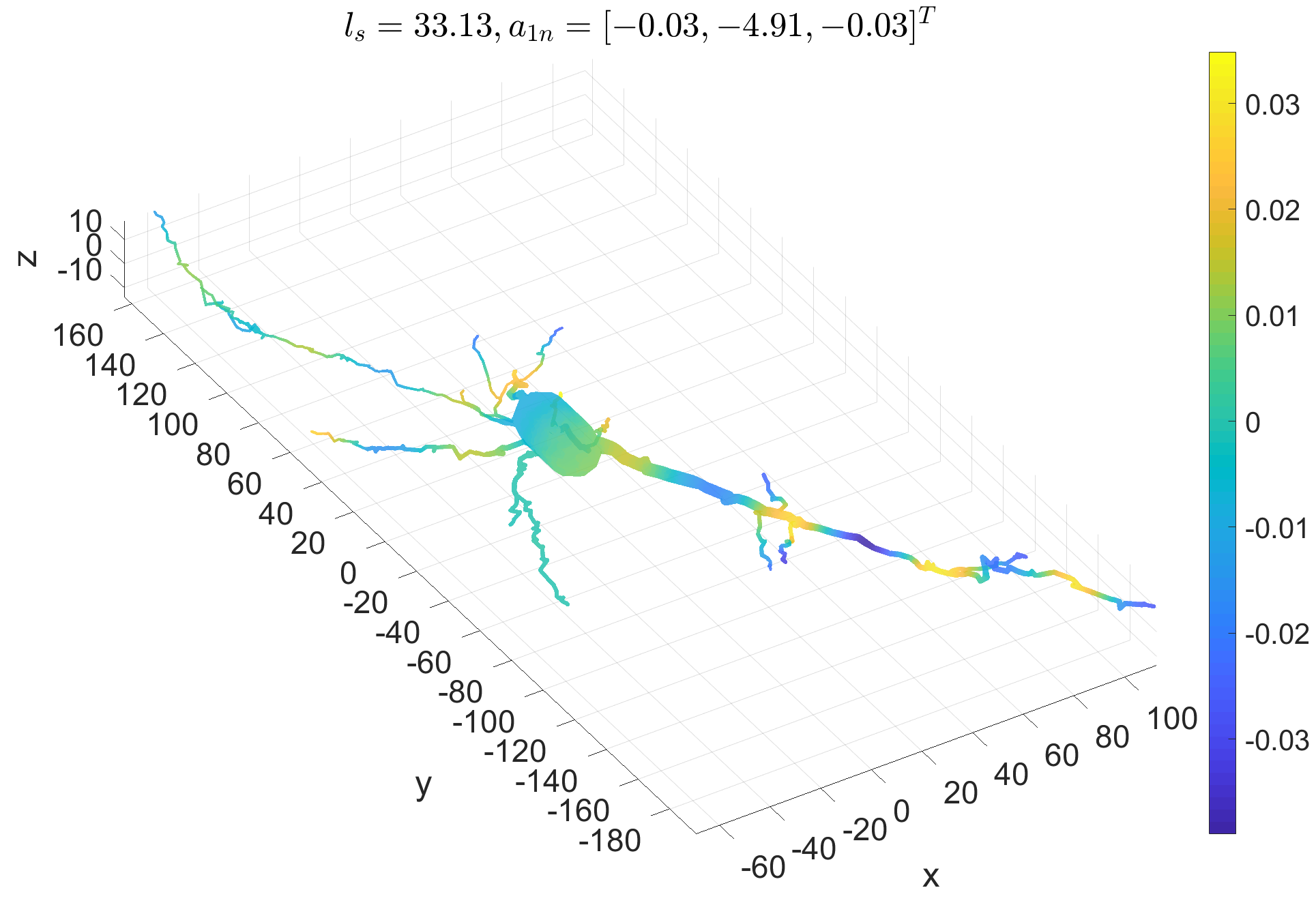}
  \caption{\label{fig:eig_big}}
\end{subfigure}~\quad
 \begin{subfigure}{0.48\textwidth}
\includegraphics[width=1\textwidth]{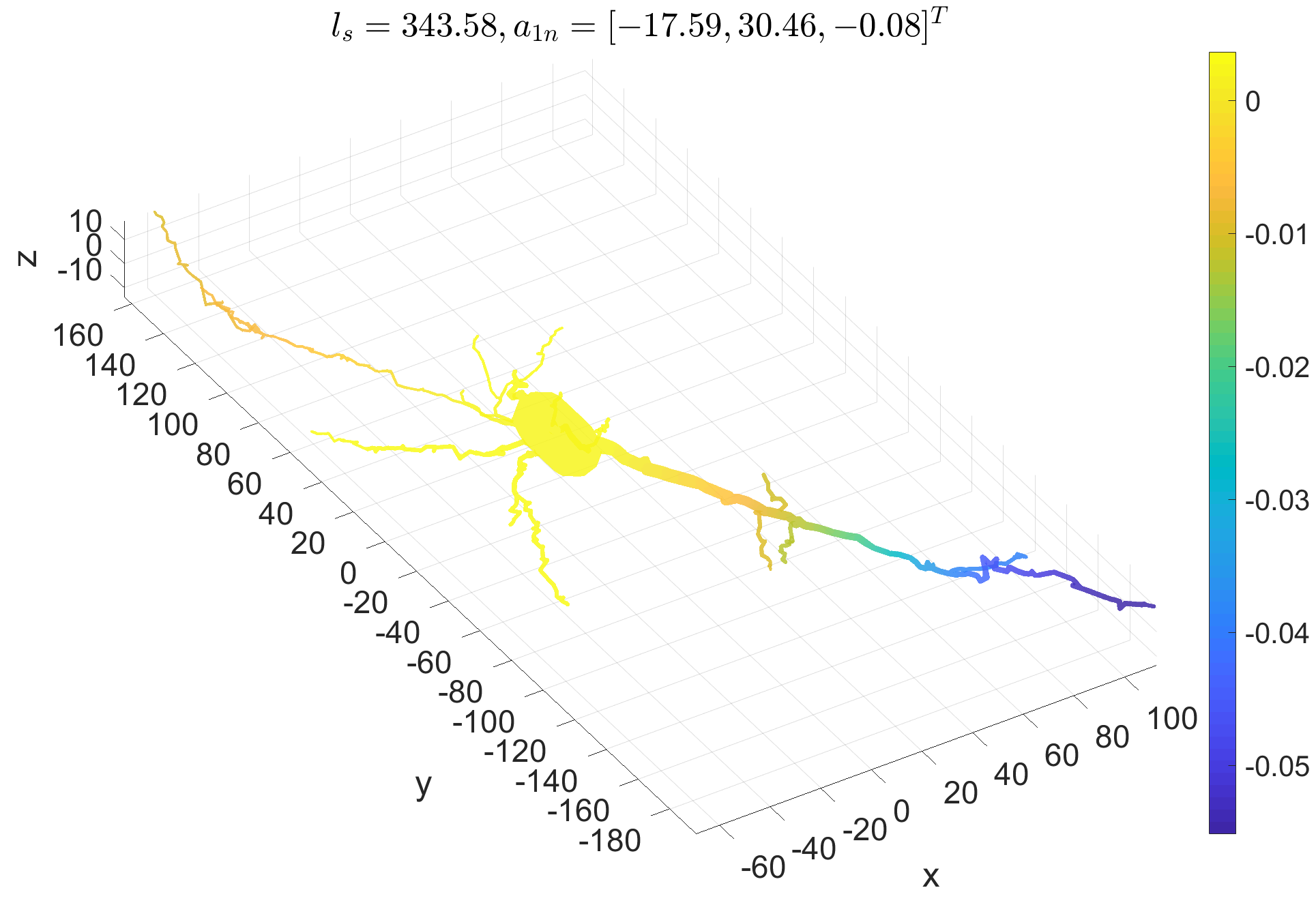}
  \caption{\label{fig:eig_long}} 
\end{subfigure}\vspace{0.4cm}
  \begin{subfigure}{0.48\textwidth}
\includegraphics[width=1\textwidth]{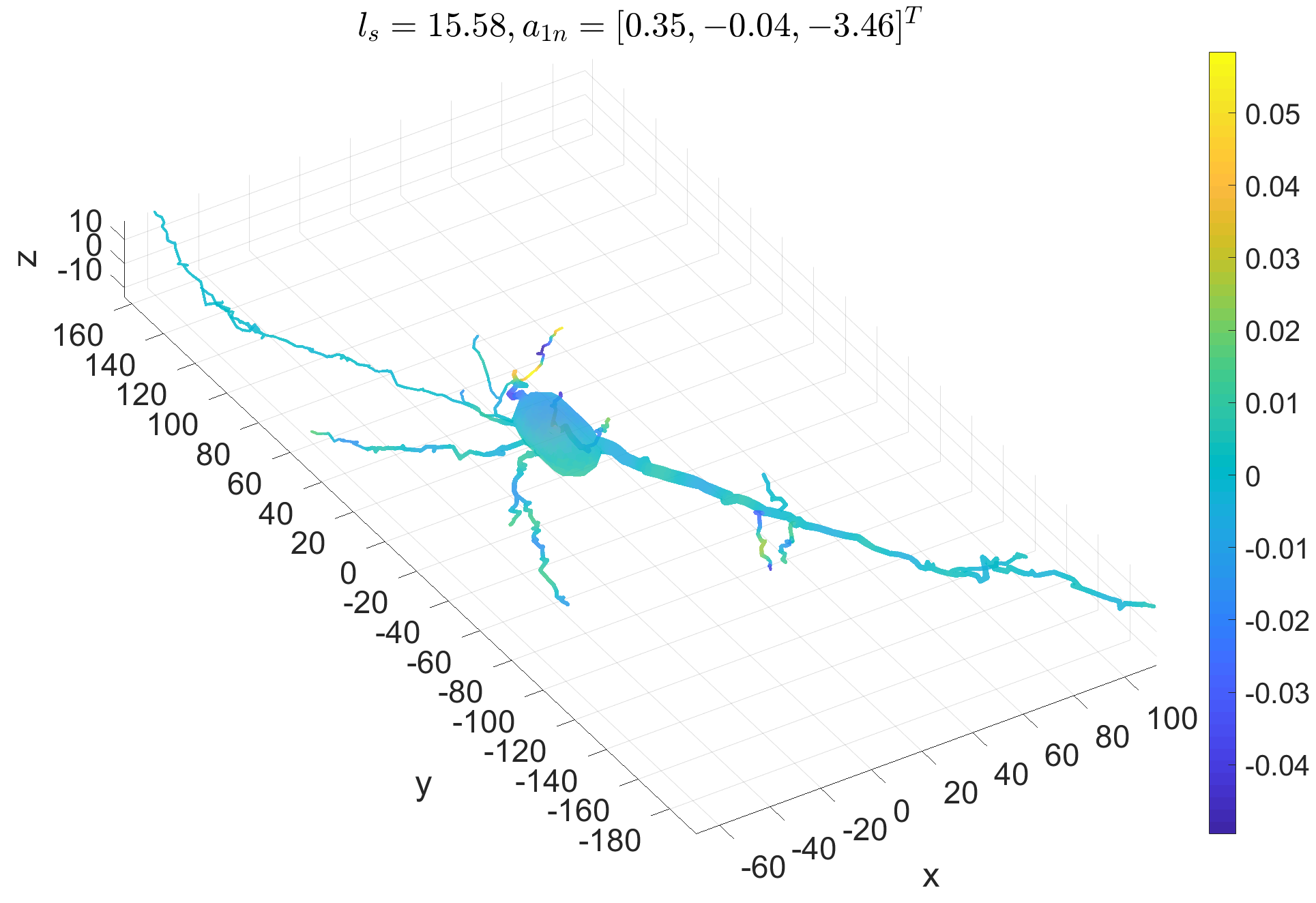}
  \caption{\label{fig:eig_z}} 
\end{subfigure}~\quad
\begin{subfigure}{0.48\textwidth}
\includegraphics[width=1\textwidth]{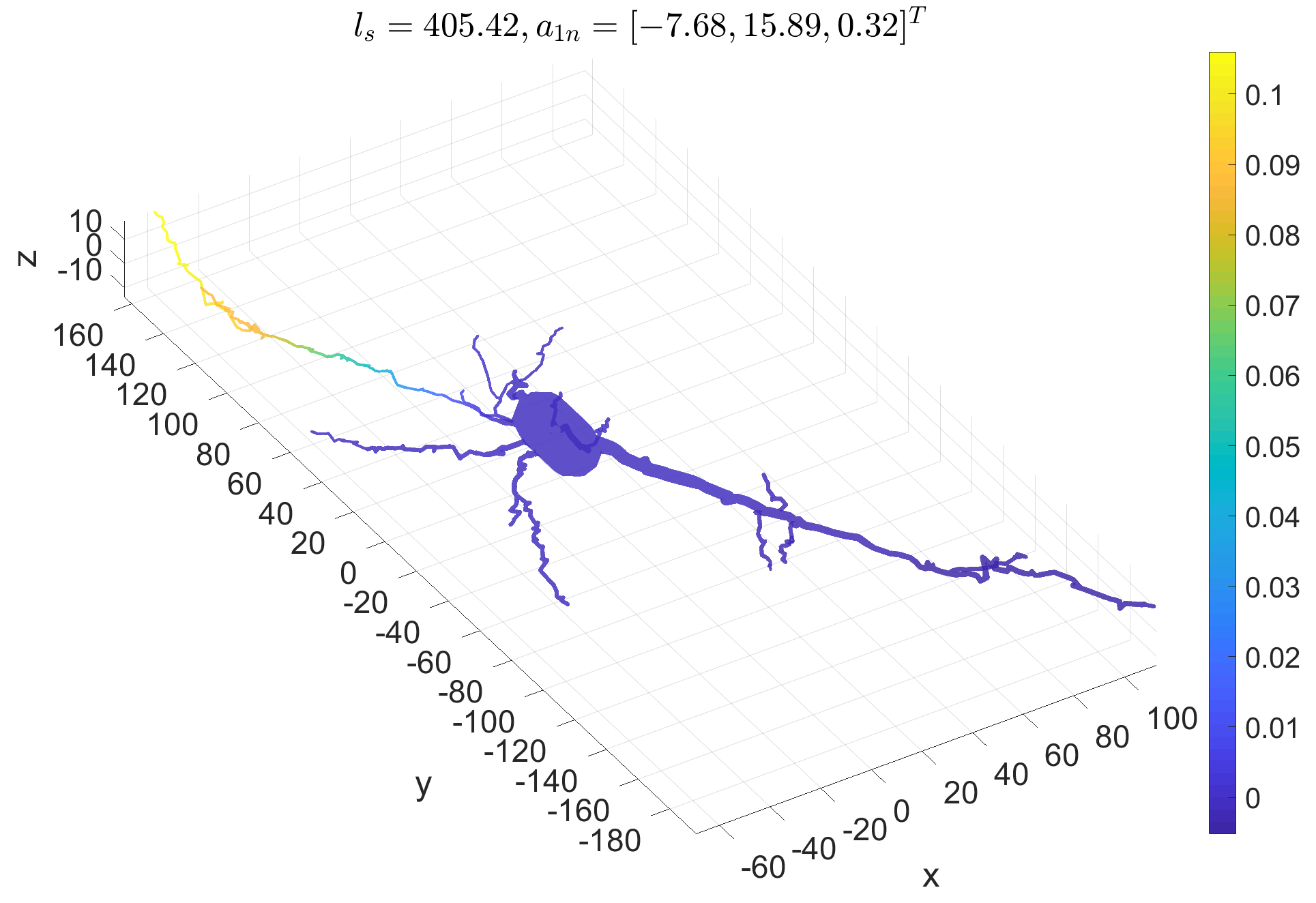}
  \caption{\label{fig:eig_longest} }
  \end{subfigure}%
\caption{ \label{fig:lapeigs} \soutnewr{}{The eigenfunctions corresponding to the spatial scale of $l_s(\lambda_n) = 33\lunit$ \textbf{(a)}, $343.6\lunit$ \textbf{(b)},  
$15.6\lunit$ \textbf{(c)} and $405\lunit$ \textbf{(d)}.  
The geometry is the pyramidal neuron {\it 02b\_pyramidal1aACC}.}}
\end{figure}

In Figure \ref{fig:sigdiff}, we show the significant eigenmodes.  
To visualize the "diffusion direction" of the eigenmodes, 
we use a RGB (red, green, blue) color scale based on the values of the
RGB vector $\bbd{c}_n$ with three non-negative valued components:
\be{}
\bbd{c}_{n}=\frac{\begin{bmatrix}
|a_{1n}^x|,  |a_{1n}^y|, |a_{1n}^z|
\end{bmatrix}}{\sqrt{(a_{1n}^x)^2+ (a_{1n}^y)^2+(a_{1n}^z)^2}}.
\ee
The color indicated by the RGB vector $\bbd{c}_n$ can be used to gauge the relative contribution of the eigenmode to the 3 principle diffusion directions, $x$, $y$, $z$.
For SEQ1, the significant eigenmodes between $25\lunit \leq l_s \leq 50\lunit$ are mostly green, meaning they contribute to diffusion in the $y$ direction.  Between $0 \leq l_s \leq 25\lunit$, there are many more significant eigenmodes 
that are red, meaning they contribute principally to diffusion in the $x$ direction.
There is only one mode that is blue, meaning it contributes significantly to diffusion in the 
$z$ direction.  This is expected since this neuron lies principally in the $x-y$ plane. 
We see also that at the higher b-value, there are more significant eigenmodes
and the signal differences are also higher, 
compared with the lower b-value. 
The most significant eigenmode is the one with length scale $l_s(\lambda_n) = 33\lunit$, it is most aligned to the $y$-direction.  
Figure \ref{fig:eig_big} shows this eigenfunction and we see the length scale corresponds to the "wavelength" of the 
significant oscillations of the eigenfunction in the geometry.
Among eigenmodes with longer length scales, $ l_s > 50 \lunit$, there are a lot fewer
significant eigenmodes than between $0 \leq  l_s \leq 50\lunit $ 
and they are mostly in the $y$ direction (being mostly green).
The eigenmode corresponding to $l_s(\lambda_n) = 343.6\lunit$ is shown in Figure \ref{fig:eig_long}.
Its removal will result in signal differences of 
\{19.1\%, 83.8\%, 8.8\%, 40.8\%\}.  The "wavelength" of this mode can be seen to be 
longer (slower oscillation) than the mode with $l_s(\lambda_n) = 33\lunit$.
For SEQ2, at the lower length scales, around $l_s = 15\lunit$, many of the significant eigenmodes for SEQ1 are no longer significant.
In Figure \ref{fig:eig_z}, we show the eigenmode
that is blue which is significant for both SEQ1 and SEQ2, corresponding to $l_s = 15.6\lunit$.
This eigenmode is predominately in the $z$-direction, and the rapid  oscillations are found in the dendrite branches.
The eigenfunction which has the longest spatial scale 
of $l_s = 405\lunit$, shown in Figure \ref{fig:eig_longest}, results in the following signal differences when it is removed: $E^{\text{RM}, i}=\{1.06\%, 6.58\%, 0.57\%, 1.22\%\}$.

\begin{figure}[ht!]
 \centering
\captionsetup[subfigure]{size=normalsize}

 \begin{subfigure}{1\textwidth}
\includegraphics[width=0.48\textwidth]{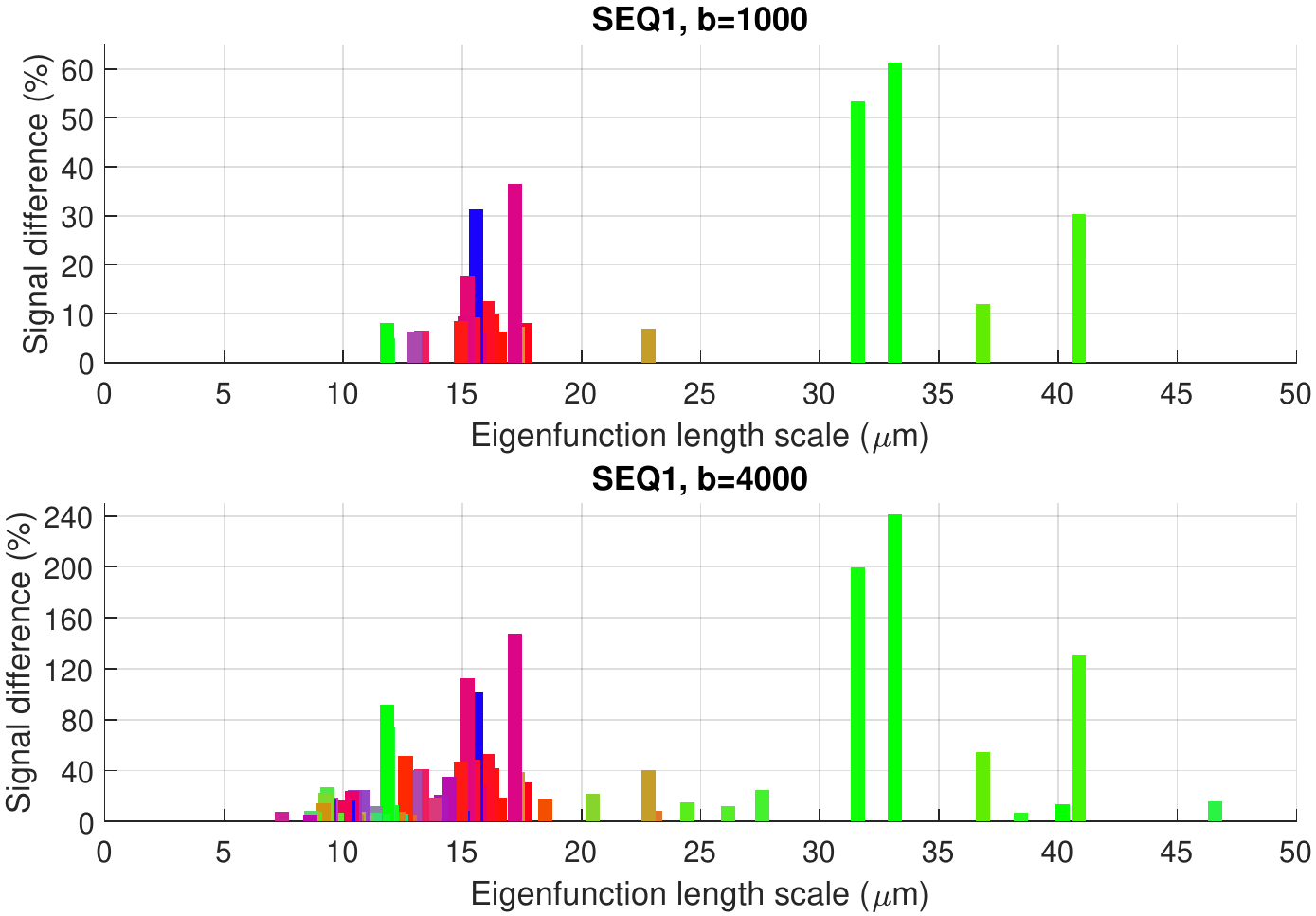} \ \
\includegraphics[width=0.48\textwidth]{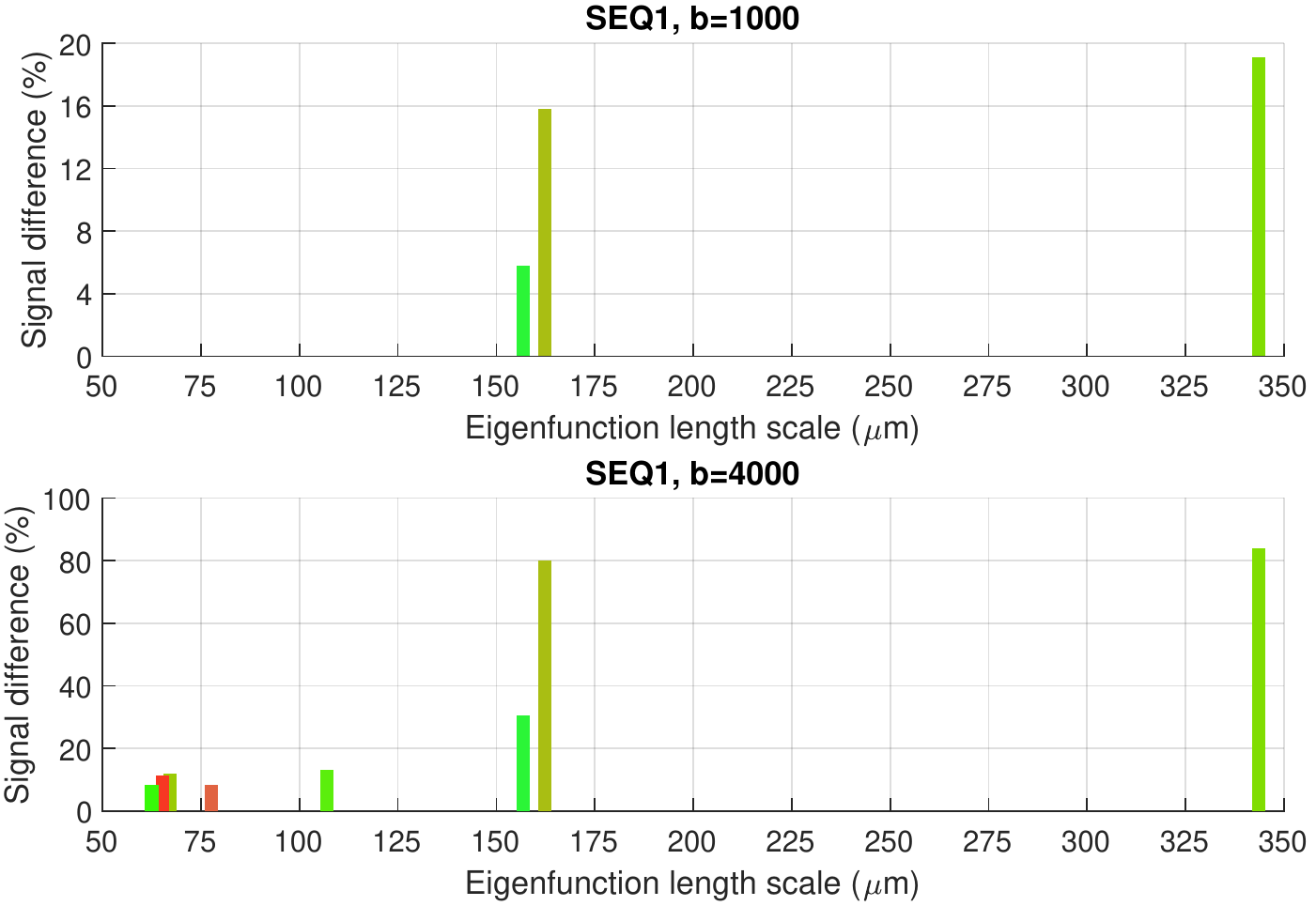}
\caption{\label{fig:sigdiff_experi1} 
SEQ1 is (PGSE, $\delta=10.6\tunit, \Delta=13\tunit$);}
\end{subfigure}%
\vspace{0.5 cm}
 
 \begin{subfigure}{1\textwidth}
\includegraphics[width=0.48\textwidth]{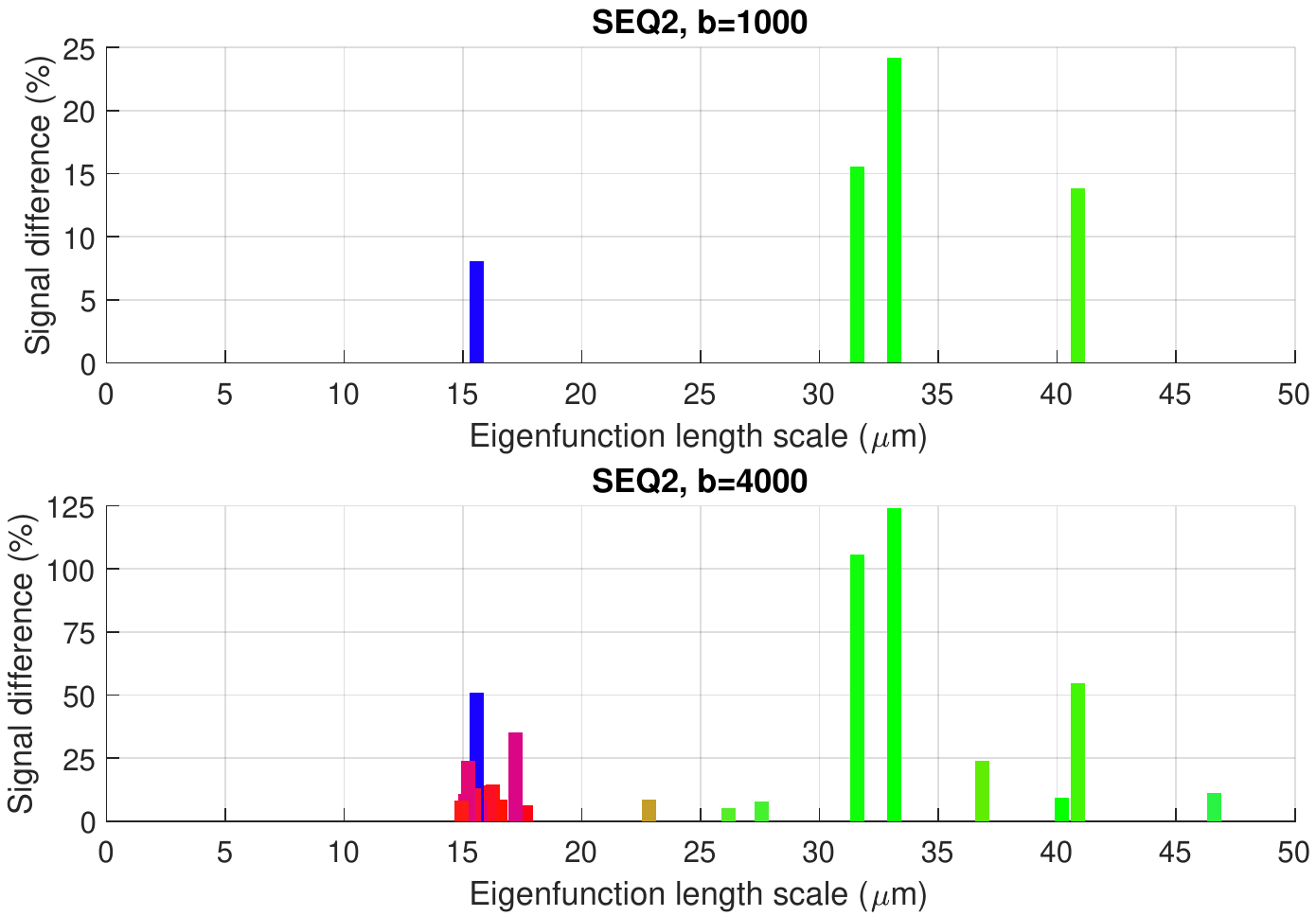} \ \
\includegraphics[width=0.48\textwidth]{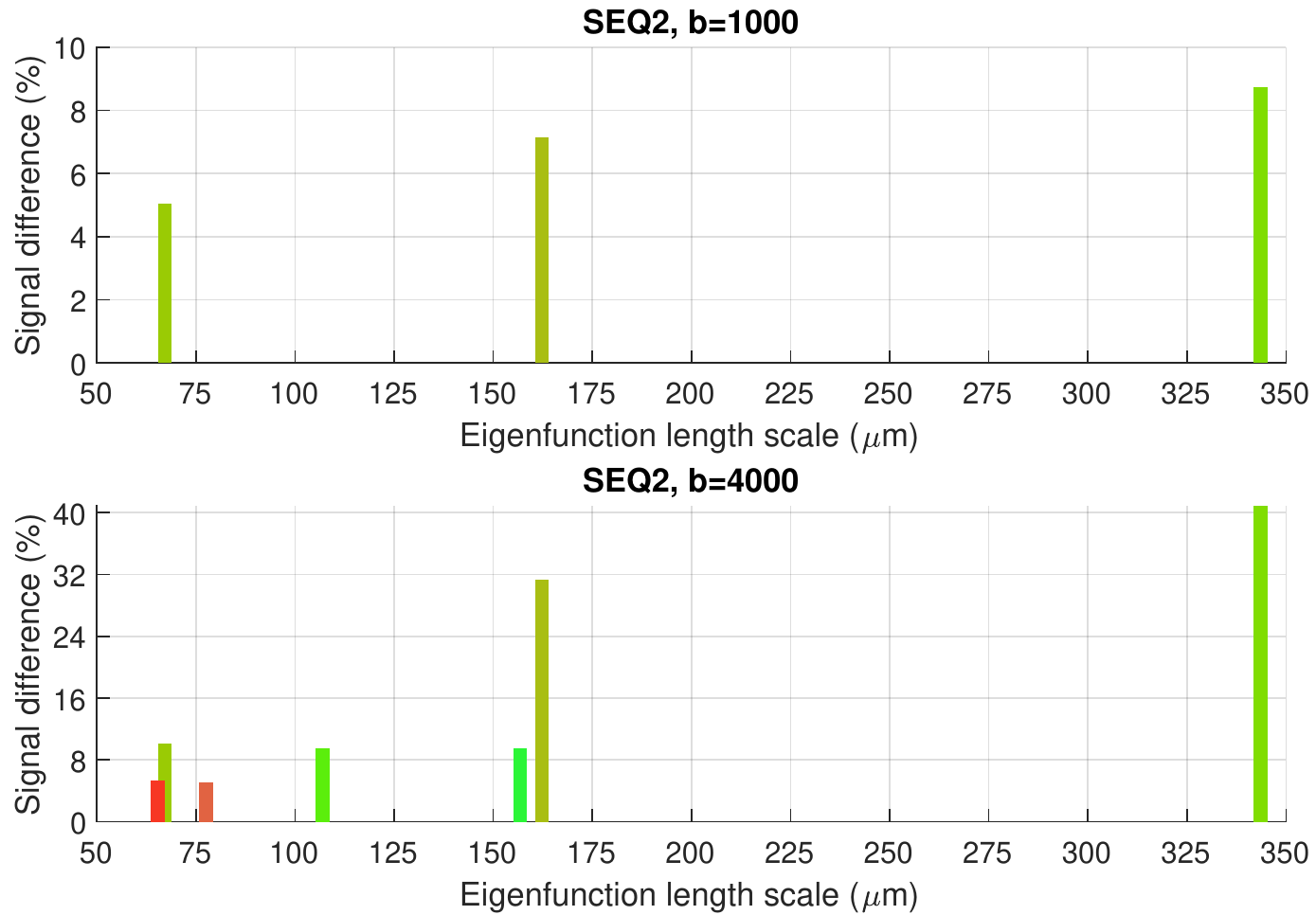}
\caption{\label{fig:sigdiff_experi3} 
SEQ2 is (PGSE, $\delta=10.6\tunit, \Delta=73\tunit$);}
\end{subfigure}

\caption{ \label{fig:sigdiff} 
The signal differences due to the removal of each eigenmode, compared to 
using the full set of 336 eigenmodes.  
The eigenvalues have been converted to a length scale
$l_s(\lambda_n)$.
The color indicates the "diffusion direction" of the eigenmodes, based on the values of the RGB vector $\bbd{c}_n$ which is related to  $\ba_{1n}$.
The geometry is the pyramidal neuron {\it 02b\_pyramidal1aACC}.
  }
\end{figure}

In Table \ref{table:mf} we summarize the number of "significant" modes 
given the threshold of $E^{\text{RM}, i}\geq 0.1\%$ and $E^{\text{RM}, i}\geq 1\%$.
Then we computed the Matrix Formalism signal using only the "significant" modes
and compared it to the reference signal  $S^\text{BTPDE}$ over 30 gradient directions.
The number of significant modes range from 27 to 197, the signal errors compared to the 
reference signal range from less than 2\% to 12\%.

\begin{table}[ht!]
\begin{center}
\begin{tabular}{|l|c|c|c|c|}
\hline 
 & Num of  &  Sig error  &Num of & Sig error   \\ 
 &   significant & from Ref&   significant  &  from Ref\\ 
 &   modes (a) & (a) & modes (b) & (b) \\ 
\hline 
SEQ1, $b=1000$ &123 & 0.5\%& 53& 12\% \\ 
\hline 
SEQ1, $b=4000$ &197 & 1.0\% & 146  &  12\% \\ 
\hline 
SEQ2, $b=1000$ &55 & 0.6\% & 27 &   9\% \\ 
\hline 
SEQ2, $b=4000$  &107 & 2.1\%& 58& 5\%\\ 
\hline 
\end{tabular} 
\caption{Significant modes are those whose removal leads to a signal difference of more than 
$0.1\%$ (a) or $1\%$ (b) compared to the signal from using the entire set of computed modes.  
In total, 336 eigenfunctions were computed. 
Signal error is the difference between the MF signal obtained using the
indicated number of significant modes compared to 
the reference signal obtained from solving the Bloch-Torrey PDE.  The signal difference 
is averaged over 30 diffusion-encoding directions, uniformly distributed on the sphere.
The geometry is the pyramidal neuron {\it 02b\_pyramidal1aACC}.
\label{table:mf}}
\end{center}
\end{table}

In Table \ref{table:timing} we give the computational times. All the simulations were performed on a server computer with 12 processors (Intel (R) Xeon (R) E5-2667 @2.90 GHz), 192 GB of RAM, running CentOS 7, using MATLAB R2019a.  It is clear that once
the eigendecomposition has been computed, the Matrix Formalism signal representation can be obtained rapidly for many sequences, b-values, and 
gradient directions. 
We note that given $n$ eigenfunctions, the number of 
associated model parameters of the Matrix Formalism representation 
is $n + 3n(n-1)/2$, because the matrix $L$ is diagonal and the three matrices $A^{i}$, $i = x,y,z$, are symmetric.
The number of parameters in each Matrix Formalism representation is also given in Table \ref{table:timing}.

\begin{table}[ht!]
\begin{center}
\begin{tabular}{|l|c|c|c|}
\hline 
&   MF & MF& BTPDE\\  
\hline
Model size& 336 modes & 197 modes &44908 nodes \\ 
  & 61656 params & 19503 params & 171017 elem \\ 
\hline
Eigen solve  & 1095 sec  &  1095 sec  &\\   
\hline 
SEQ1, $b=1000$& 0.09 sec& 0.04 sec &  26 sec\\ 
\hline 
SEQ1, $b=4000$& 0.09 sec&  0.04 sec &39 sec\\ 
\hline 
SEQ2, $b=1000$&   0.12 sec & 0.05 sec & 22 sec\\ 
\hline 
SEQ2, $b=4000$&  0.13 sec  & 0.05 sec &30 sec\\ 
\hline 
\end{tabular} 
\caption{The computational times to calculate the signals due to the indicated diffusion-encoding 
sequences at the indicated b-values, averaged over 30 diffusion-encoding directions.
The geometry is the pyramidal neuron {\it 02b\_pyramidal1aACC}.
\label{table:timing}} 
\end{center}
\end{table}

\subsection{Eigenvalues and eigenfunctions of the Bloch-Torrey operator}
%\marginparnew{This Section is new}
At high gradient amplitudes, it was demonstrated\cite{Grebenkov2014,Moutal2019} in certain geometries (intervals, disks, spheres and the exterior of arrays of disks), the magnetization solution of the Bloch-Torrey equation exhibits localization near boundaries and interfaces.
The analysis in those papers was based on the eigenfunctions of the complex-valued Bloch-Torrey operator, 
\ben
-\left(\Dintr \nabla ^2 - \bi\gamma  \bg \cdot \bx\right),
\een
in contrast to the Laplace operator $-\Dintr \nabla ^2 $.
In the Appendix \ref{sec:appendix2}, we show that the conversion between the eigenfunctions of the Bloch-Torrey operator $\{\psi_j\}$ and 
the eigenfunctions of the Laplace operator $\{\phi_j\}$ is given by
\be{}
\begin{split}
[ \psi_1(\bx),\cdots\psi_{N_{eig}}(\bx)]^T &= V^{-1} [ \phi_1(\bx),\cdots\phi_{N_{eig}}(\bx)]^T \\
V [ \psi_1(\bx),\cdots\psi_{N_{eig}}(\bx)]^T &= [ \phi_1(\bx),\cdots\phi_{N_{eig}}(\bx)]^T 
\end{split},
\ee
where the columns of $V$ contains the eigenvectors of the complex-valued matrix $K(\bg)$.  
The eigenvalues of the Bloch-Torrey operator are exactly the eigenvalues of $K(\bg)$.  
We remind the reader that the eigenvalues of 
$K(\bg)$,
\ben 
\mu_1,\cdots,\mu_{N_{eig}}, \quad 0 <  \Re{\mu_1} \leq \cdots \leq  \Re{\mu_{N_{eig}}},
\een 
are found on the diagonal of the matrix $\Sigma$.  We order the Bloch-Torrey (BT) eigenvalues and eigenfunctions
by the magnitude of the real part of $\mu_{i}$, all of which are strictly greater than $0$.

In this section, we simulate the spindle neuron {\it 03b\_spindle4aACC}.  
We compute $\psi_1,\cdots,\psi_{N_{eig}}$ using the transformation above.
To better visualize the "localized" nature of the BT eigenfunctions, in the Figures, 
we only indicate the 
support of $\psi_i(\bx)$:
\ben
\text{supp}_{\psi,\epsilon} \equiv  \{\bx: \vert \psi(\bx)\vert \geq 0.01 \max_{\bx} \vert \psi\vert \},
\een
i.e., the region in space where it attains at least 1\% of its maximum value.  
We do not indicate the 
actual values of $\psi$ inside its support.

\begin{figure}[ht!]
 \centering
\captionsetup[subfigure]{size=normalsize}
 \begin{subfigure}{0.95\textwidth}
\includegraphics[width=0.47\textwidth]{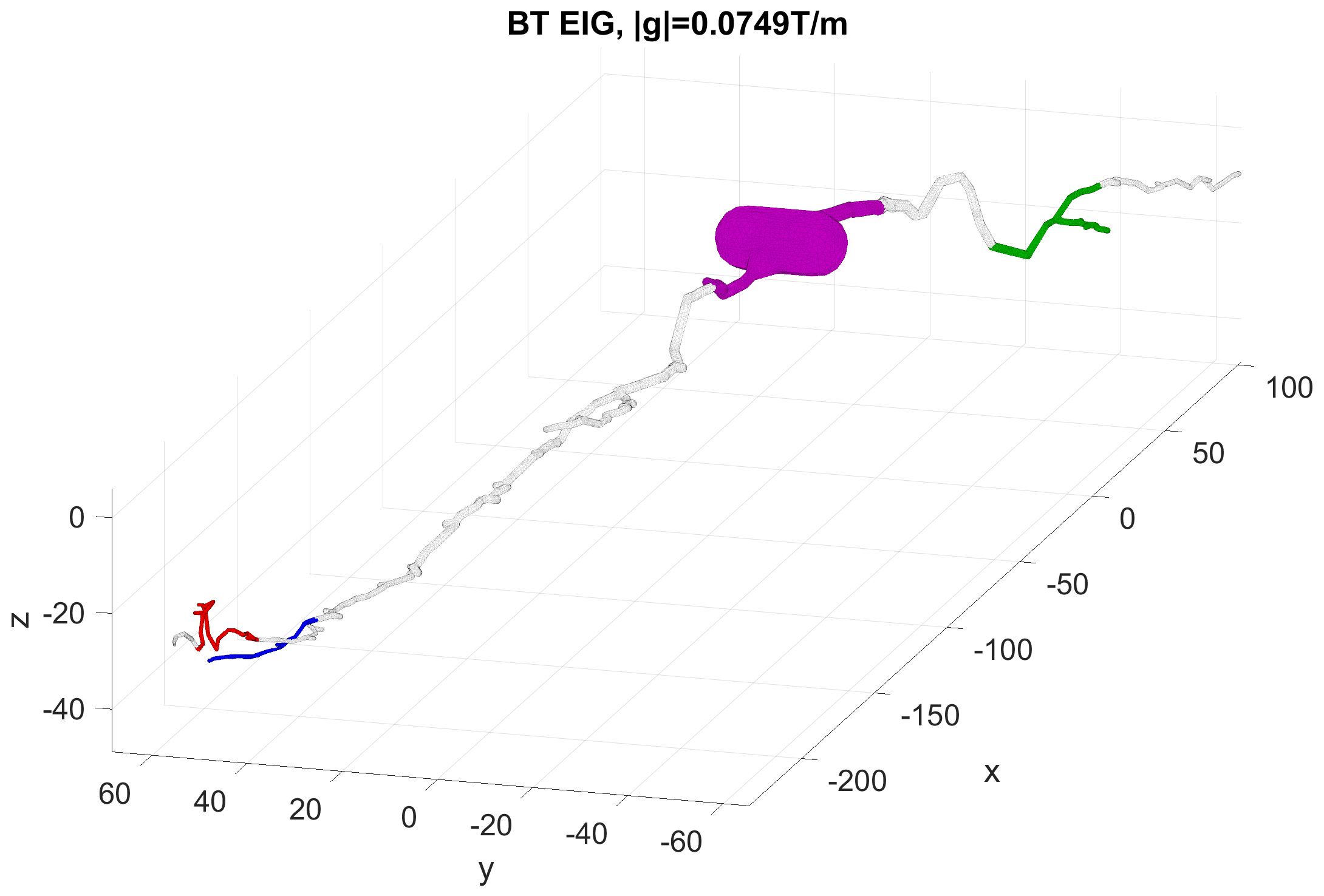}
\includegraphics[width=0.52\textwidth]{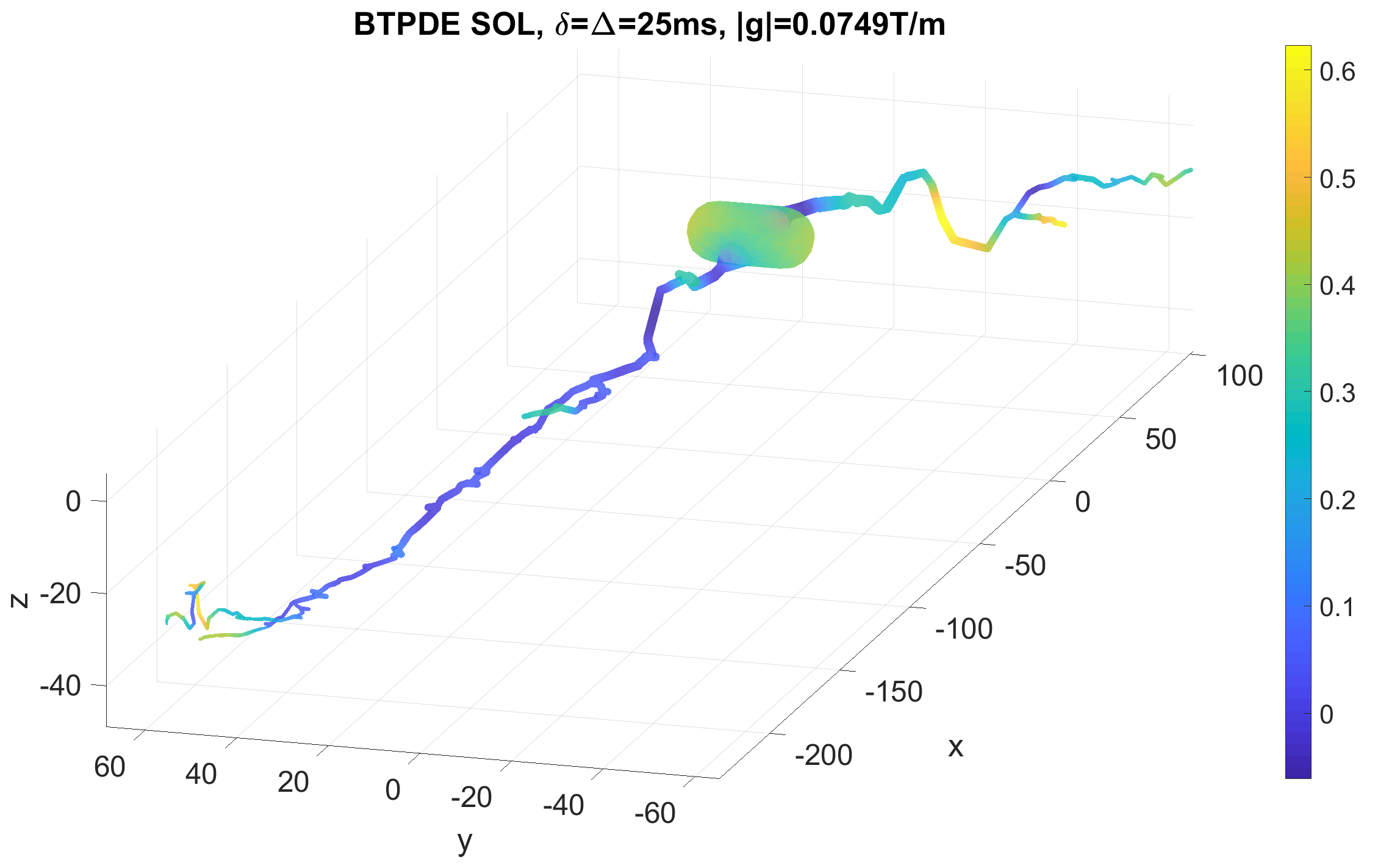}
\caption{\label{fig:eig_2e5} $\vert \bg \vert = 0.075 \qunit$}
\end{subfigure}%
\vspace{0.3 cm}
 \begin{subfigure}{0.95\textwidth}
\includegraphics[width=0.47\textwidth]{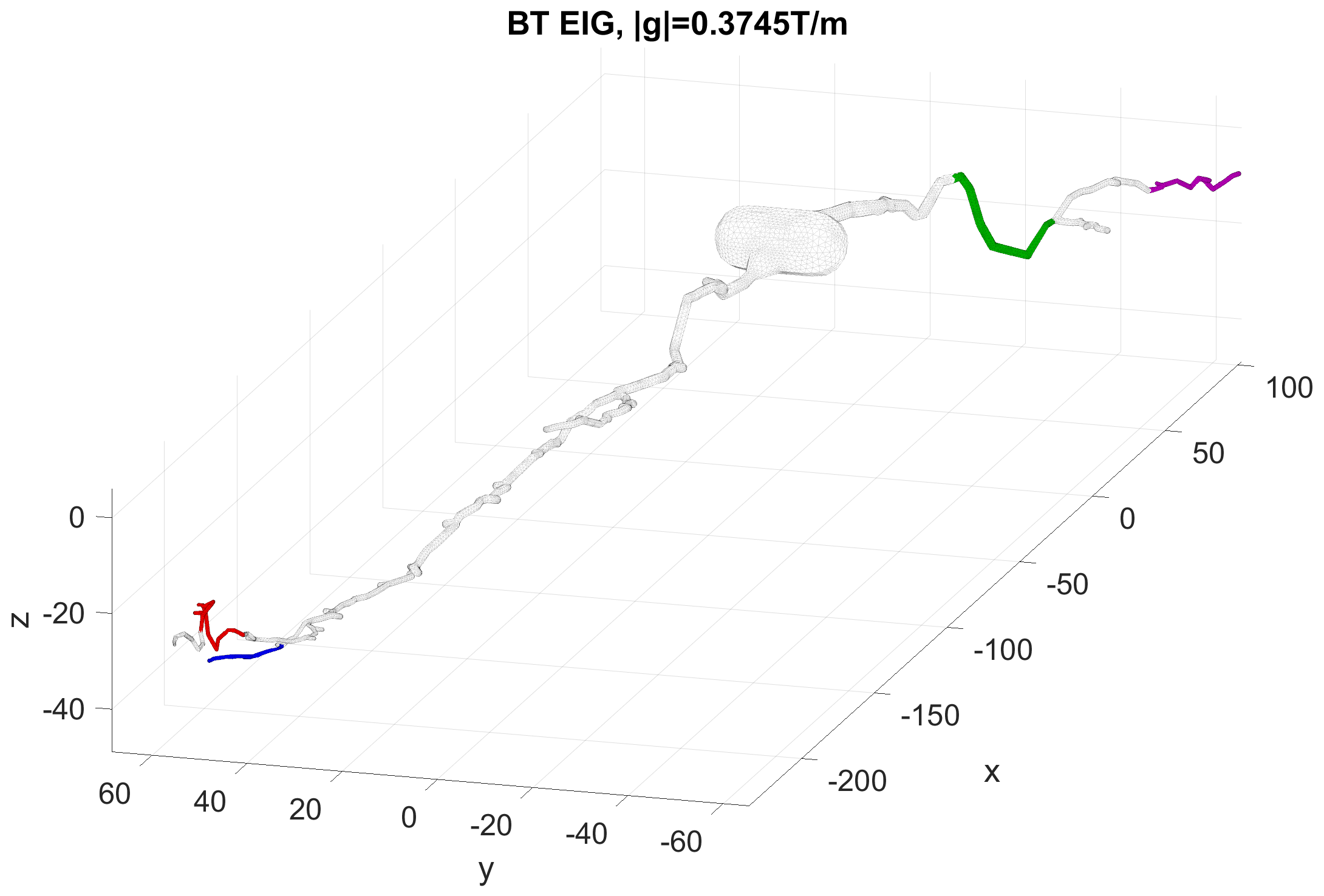}
\includegraphics[width=0.52\textwidth]{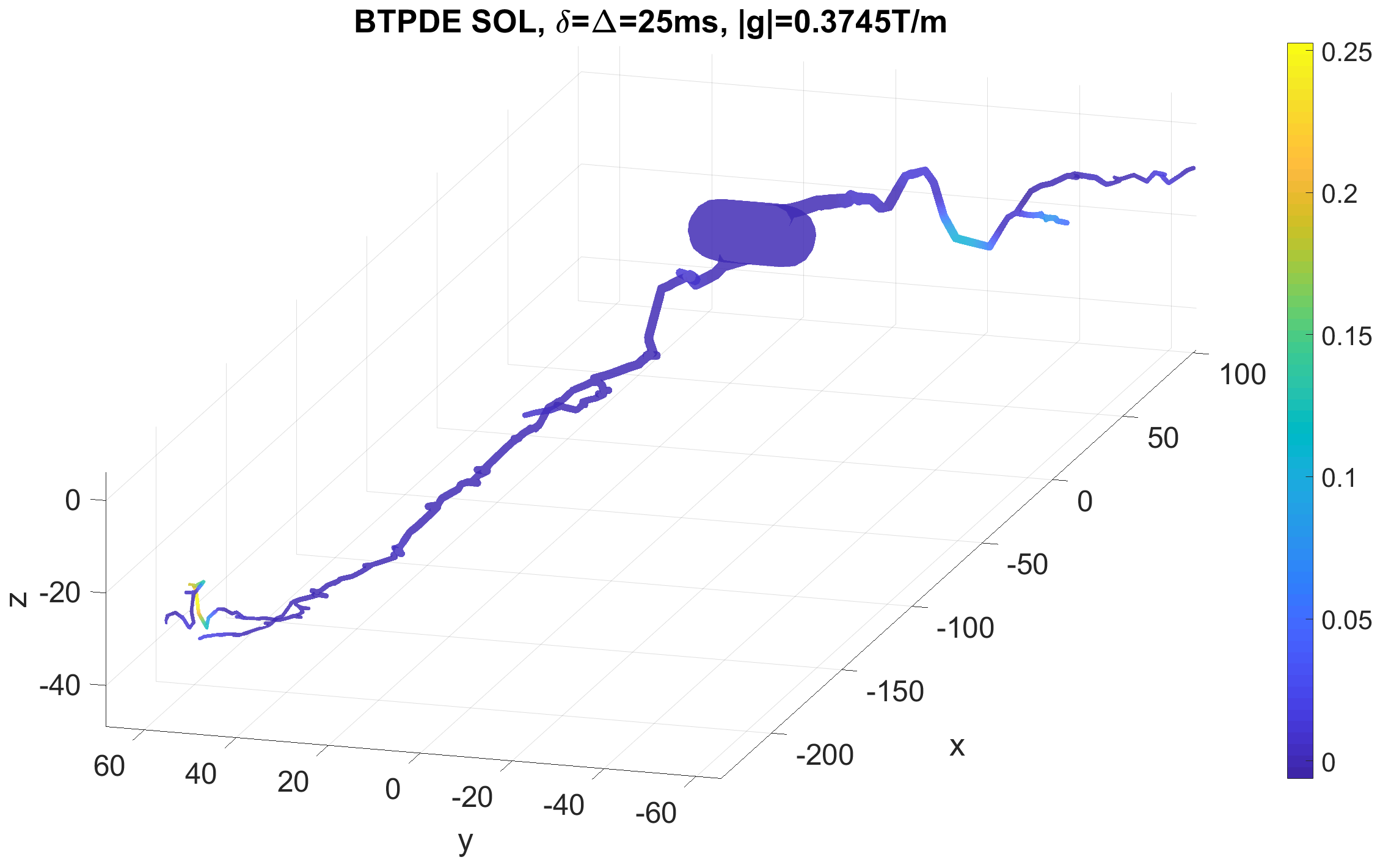}
\caption{\label{fig:eig_1e4}$\vert \bg \vert = 0.3745 \qunit$}
\end{subfigure}
 \caption{ \label{fig:BT_eig} 
\textbf{(a)}  The gradient amplitude is $\vert \bg \vert = 0.075 \qunit$.
Left: the supports of the first few BT eigenfunctions, ${\psi_1}$ (red), ${\psi_2}$ (green),
${\psi_4}$ (blue), ${\psi_5}$ (violet).   
Right: the magnetization at $T_E = \delta+\Delta$, where $\delta=\Delta=25\tunit$.
\textbf{(b)} The gradient amplitude is $\vert \bg \vert = 0.3745 \qunit$.
Left: the supports of the first few BT eigenfunctions, ${\psi_1}$ (red), ${\psi_2}$ (green),
${\psi_5}$ (blue), ${\psi_6}$ (violet).    
Right: the magnetization at $T_E = \delta+\Delta$, where $\delta=\Delta=25\tunit$.
The BT eigenfunction index is in order 
of the magnitude of the real part.  
The gradient direction is $\bug =[0.7071  \ \ 0  \ \  0.7071]$. 
The geometry is the neuron {\it 03b\_spindle4aACC}.
  }
\end{figure}

In Figure \ref{fig:eig_2e5}, at the gradient amplitude $\vert \bg \vert = 0.0749$T/m, we show the supports of 
the BT eigenfunctions numbered 1, 2, 4, 5.  We see 
the support of $\psi_5$ lies mostly in the soma.  The supports of 
$\psi_1$, $\psi_2$, $\psi_4$ lie on the dendrite branches.  The support of $\psi_3$ (lying on the dendrite branches) is not shown to 
avoid overlapping support regions in the plot.  
The non-zero values of the magnetization at $T_E$ show the locations of the supports of BT eigenfunctions 1, 2, 4, 5.
\soutnewr{}{In Figure \ref{fig:eig_1e4},} at the gradient amplitude $\vert \bg \vert = 0.3745$T/m, the supports of the BT eigenfunctions 
numbered 1, 2, 5, 6 all lie on the dendrites.  The supports of $\psi_3$ and $\psi_4$ (lying on the dendrite branches) 
are not shown to avoid overlapping regions in the plot.  
The non-zero values of the magnetization show the locations of the supports of these eigenfunctions.  

\begin{figure}[ht!]
 \centering
\includegraphics[width=0.48\textwidth]{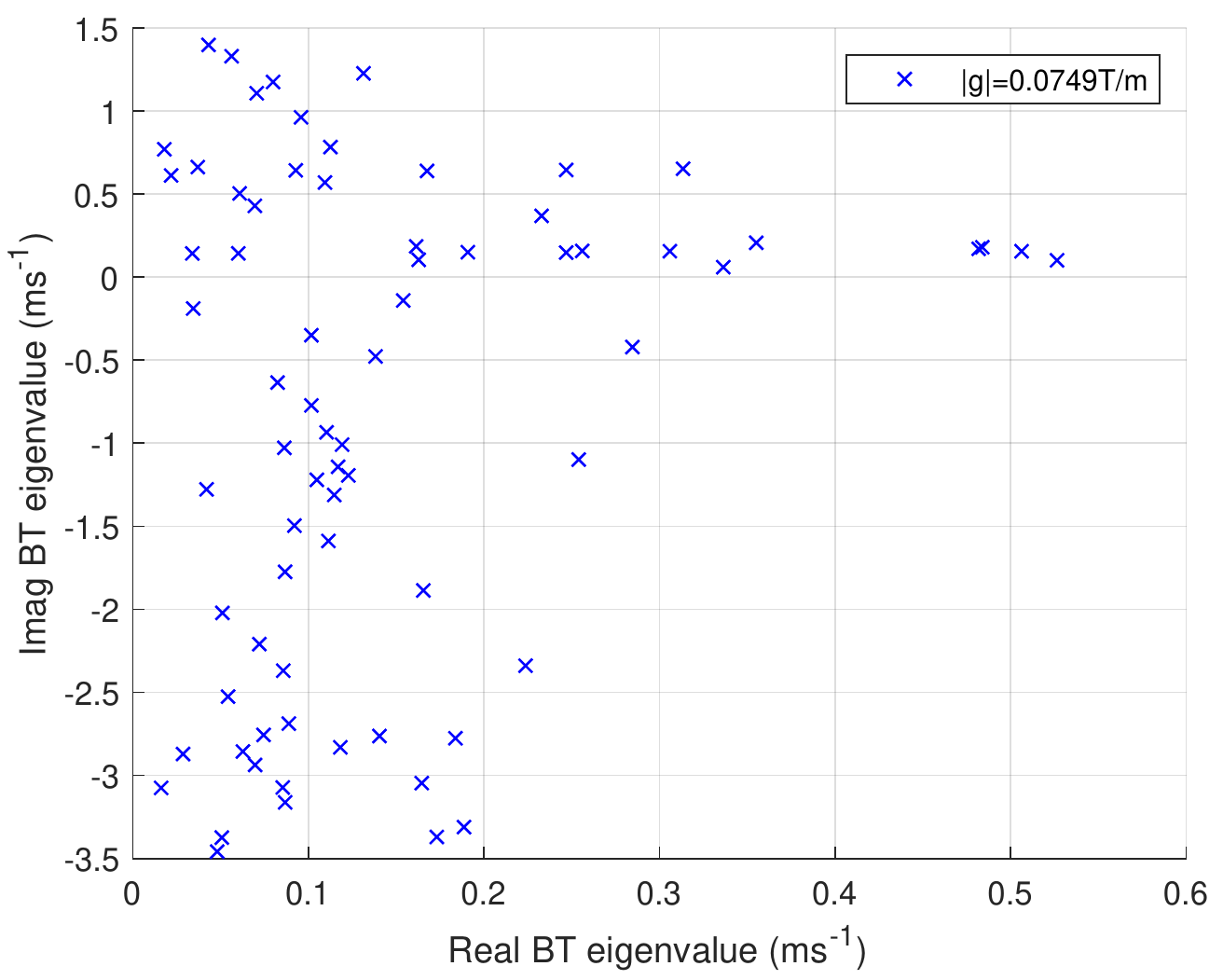}
\includegraphics[width=0.48\textwidth]{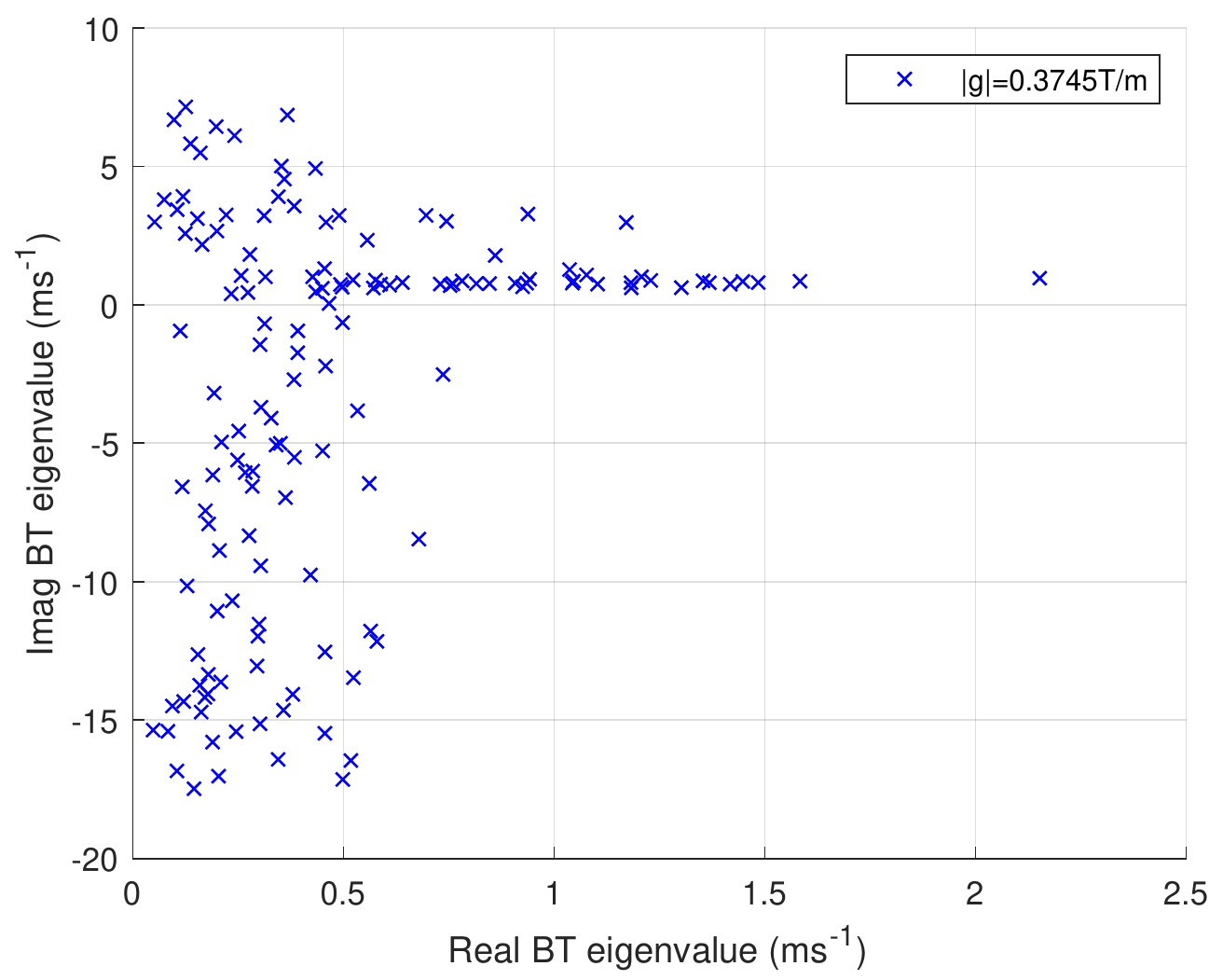}
  \caption{The significant BT eigenvalues after the projection of the
	initial spin density onto the space of BT eigenfunctions.  The real part of the BT eigenvalues 
are plotted against the imaginary part.  
	The gradient direction is $\bug =[0.7071  \ \ 0  \ \  0.7071]$.
	\textbf{Left:} $\vert \bg \vert = 0.075 \qunit$.  \textbf{Right:} $\vert \bg \vert = 0.3745 \qunit$.
	The geometry is the neuron {\it 03b\_spindle4aACC}.
	\label{fig:btcoeff_t0}}
\end{figure}

Next, we projected the initial spin density, $\rho$, which is the constant function, onto the BT eigenfunctions:
\ben
\rho \sqrt{\vert \Omega \vert }\left( V_{11}\psi_1+\cdots+V_{1N_{eig}}\psi_{N_{eig}}\right).
\een
The projection coefficients are proportional to the first row of the matrix $V$.  
We normalized the BT eigenfunctions so that they have unit $L^2$ norm 
(integral of the square of the function over the geometry).
We call a BT eigenfunction $\psi_{j}$ {\it significant} if 
\ben
\vert V_{1j}\vert \geq 0.01.
\een
In Figure \ref{fig:btcoeff_t0} we show the complex eigenvalues of the {\it significant} BT
eigenfunctions.  We see that at the higher gradient amplitude, the
significant BT eigenvalues have a wider range in both their real parts as well as their imaginary parts, 
than at the smaller gradient amplitude.  In terms of the magnetization, a larger range of the real part indicates faster
transient dynamics during the gradient pulses, and a larger range of the imaginary part indicates more 
time oscillations.  

\subsection{Length scale range of Laplace eigenfunctions}
%\marginparnew{This Section is completely new}

To investigate the question of what range of Laplace eigenvalues is sufficient to accurately describe 
the diffusion MRI signal, we computed $S^{\text{BTPDE}}$ and 
$S^{\text{MF}}$ for the neuron {\it 03b\_spindle4aACC}
on two finite elements meshes, a fine mesh with 44201 nodes and 160205 elements, and a coarser mesh 
with 17370 nodes and 56163 elements.  The $S^{\text{MF}}$ on these two FE meshes were computed 
with 4 choices of $l_s^{min}$: $4 \lunit$, $2 \lunit$, $1 \lunit$ and $0.5 \lunit$.  The $S^{\text{BTPDE}}$ on the finer mesh is considered the reference solution and we show in Table \ref{table:eiglength_range} the mean
relative signal error between the various simulations and the reference solution, averaged over the 6 
gradient directions.  Two gradient amplitudes and three PGSE sequences (with $\delta=\Delta$, $\delta = 2.5\tunit, 10\tunit, 25 \tunit$) were simulated, making a range of b-values from $4.2\bunit$ to $10416\bunit$.
We see that on the coarser finite elements mesh, the mean relative error is less than $4.1\%$ for 
$S^{\text{MF}}$ with $l_s^{min}=4\lunit$ and it is less than $5.4\%$ 
for $S^{\text{MF}}$ with $l_s^{min}=0.5\lunit$ for the entire range of b-values.  
On the fine finite elements mesh, the mean relative errors are less than 
$2.6\%$, $1.3\%$, $1.6\%$, $1.6\%$ for $S^{\text{MF}}$ with 
$4 \lunit$, $2 \lunit$, $1 \lunit$ and $0.5 \lunit$, respectively.
Some small increase in the mean relative errors are probably due to numerical errors of 
the eigenvalue computations.

\begin{table}[ht!]
	\centering
	\begin{tabular}{|c|c|c|c||c|c|c||c|c|c|c|}
		\hline
\multirow{2}{*}{ }Gradient& Pulse & b-value &BTPDE-h (REF) & \multicolumn{7}{c|}{Mean relative error with respect to BTPDE-h (REF)
(\%)} \\ \cline{5-11}
 amplitude &duration  & ($\bunit$) &  Mean signal/$S_0$ & MF-H 4 &MF-H 0.5 & BTPDE-H &MF-h 4 &MF-h 2 &MF-h 1 & MF-h 0.5   \\ \hline\hline
\multirow{3}{*}{ $0.075 $ T/m} & 
$2.5$ ms & 4.2& 0.995             &0.0 & 0.0& 0.0 & 0.0 &0.0 &0.0 &0.0 \\ \cline{2-11}
& $10 $ ms& 267& 0.829        &0.4  & 0.5 &0.0 & 0.1& 0.1& 0.1& 0.1\\\cline{2-11}
& $25 $ ms & 4167 &  0.262    &2.4 & 2.6& 0.2 & 0.6 &0.7 &0.7& 0.7\\\hline \hline
\multirow{3}{*}{ $ 0.3745 $ T/m }& 
$2.5 $ ms & 104& 0.884           & 0.3 & 0.6& 0.1  &0.1& 0.1& 0.1 &0.1\\\cline{2-11}
& $10 $ ms & 6667& 0.073      &4.1 & 5.4 &0.3 & 0.9 &1.3 &1.6& 1.6\\\cline{2-11} 
& $25 $ ms & 104167& 0.009 &1.9 & 2.2& 1.1 & 2.6& 0.7 &0.9 &1.1 \\\hline 
\end{tabular}
\caption{Diffusion MRI signals were simulated in 6 gradient directions, uniformly distributed 
in the unit semi-sphere.  The reference signal, denoted BTPDE-h, is the BTPDE signal on the finely discretized 
finite element mesh (44201 nodes, 160205 elements), with ODE solver tolerances $atol = 10^{-6}$,
$rtol=10^{-4}$.  The mean signal is the averaged signal over the 6 gradient directions.  
The mean relative error is the relative difference between the various simulations 
and the reference signal BTPDE-h, averaged over the 6 gradient directions.
The simulation denoted as BTPDE-H is the BTPDE solution on a coarser finite elements mesh 
(17370 nodes, 56163 elements) with the same ODE solver tolerances.
Two Matrix Formalism simulations on the coarse FE mesh are denoted "MF-H 4" ($l_s^{min}=4\lunit$, 166 eigenfunctions)
and "MF-H 0.5" ($l_s^{min}=0.5\lunit$, 5498 eigenfunctions).
Four Matrix Formalism simulations on the fine FE mesh are denoted "MF-h 4" ($l_s^{min}=4\lunit$, 170 eigenfunctions),
 "MF-h 2" ($l_s^{min}=2\lunit$, 507 eigenfunctions), "MF-h 1" ($l_s^{min}=1\lunit$, 2082 eigenfunctions),
and "MF-h 0.5" ($l_s^{min}=0.5\lunit$, 9093 eigenfunctions).
The geometry is the neuron {\it 03b\_spindle4aACC}. 
\label{table:eiglength_range}}
\end{table}

\begin{figure}[ht!]
 \centering
\includegraphics[width=0.31\textwidth]{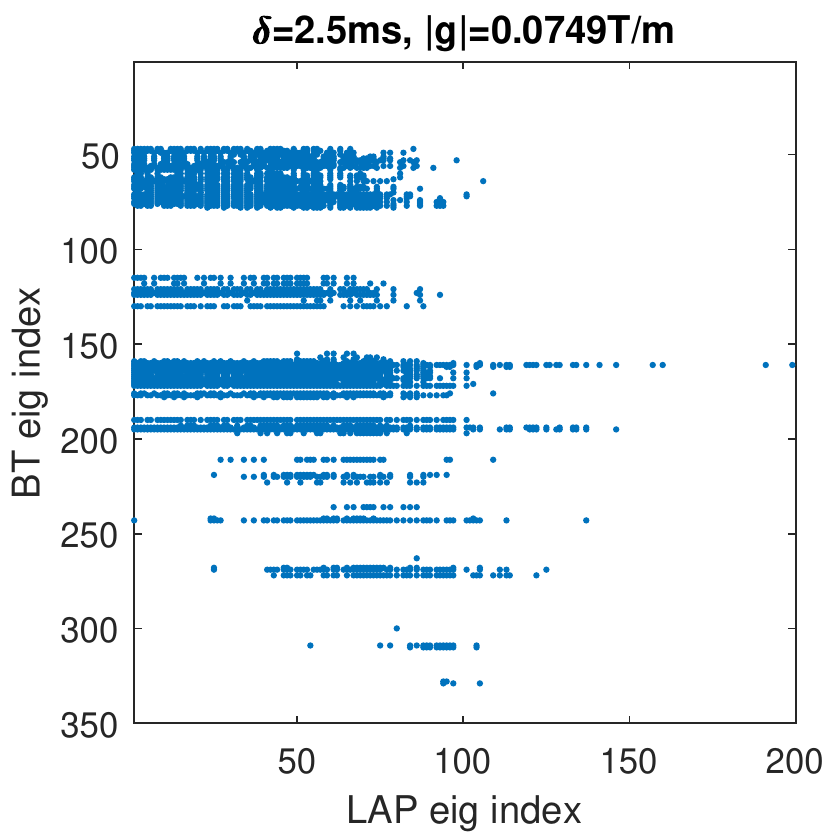} \quad
\includegraphics[width=0.31\textwidth]{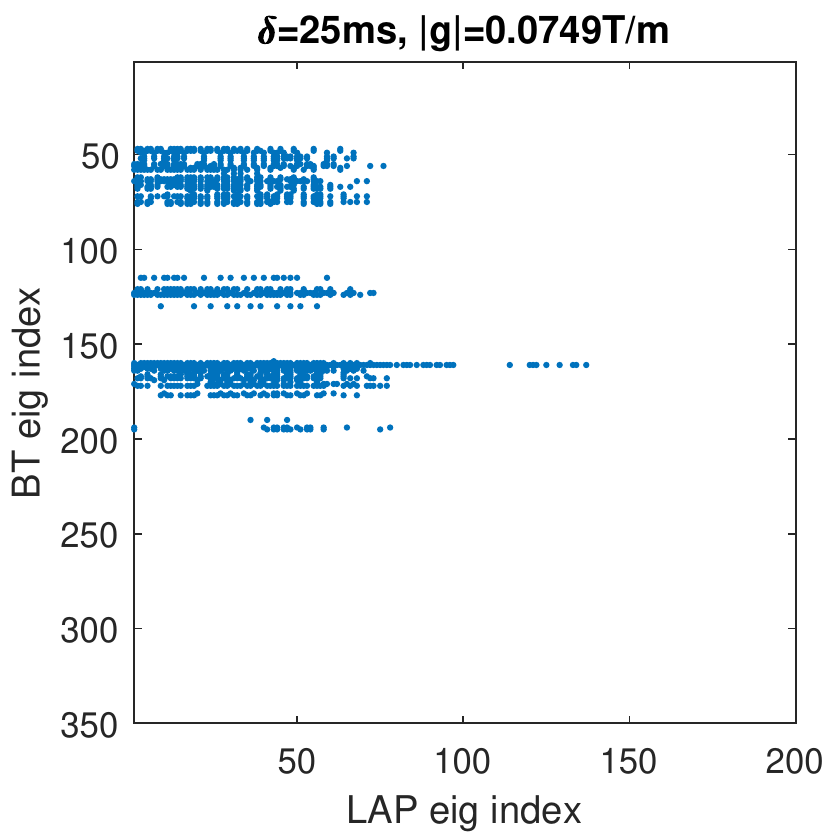} \quad
\includegraphics[width=0.31\textwidth]{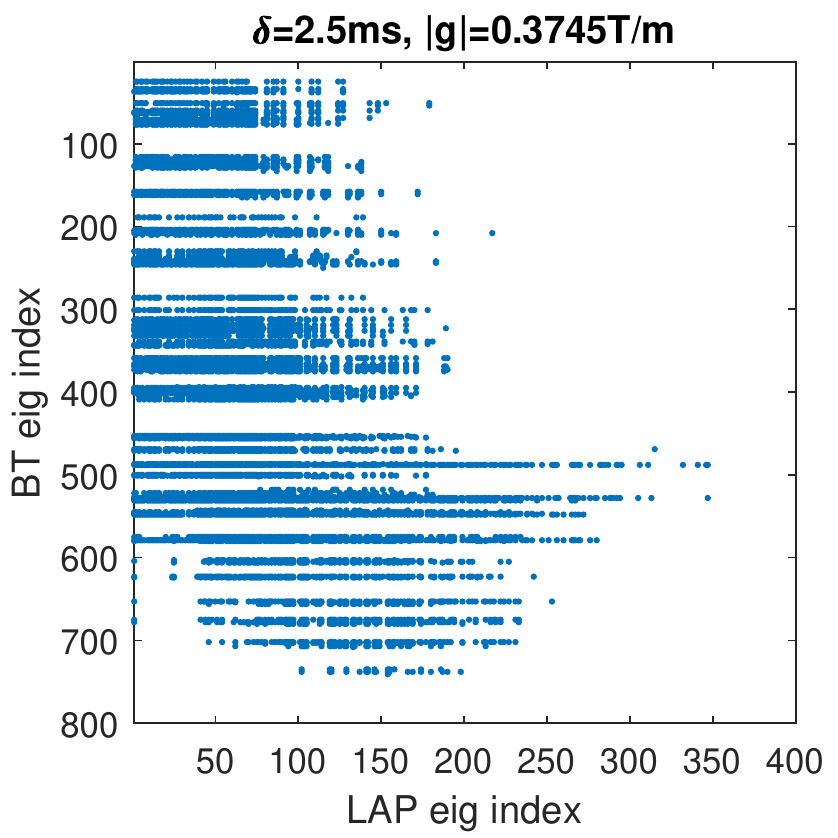}
\caption{\label{fig:lapbt_grid}
The magnetization at $t=\delta$ is projected onto the space of BT eigenfunctions, each of which 
is projected onto the Laplace eigenfunctions.  
The matrix that contains the projection coefficients of the time-evolved BT eigenfunctions onto the
Laplace eigenfunctions is shown here.
If the projection of the time-evolved $i$th BT eigenfunction onto the $j$th Laplace eigenfunction 
has a magnitude greater than 0.001, it is marked by a 
blue dot in the $i$th row and $j$ column.  
The BT and Laplace eigenfunction indices are ordered by the real part of their eigenvalues.
\textbf{Left:} $\vert \bg \vert = 0.075 \qunit$, $\delta = \Delta= 2.5\tunit$, $b=4.17\bunit$, the maximum BT index is 334, 
corresponding to the eigenvalue of $\mu = (3.3+0.1\bi)\tunit^{-1}$, the maximum LAP index is 199, 
corresponding to the length scale of $3.5\lunit$.
\textbf{Middle:} $\vert \bg \vert = 0.075 \qunit$, $\delta = \Delta = 25\tunit$, $b=4167\bunit$ the maximum BT index is 195, 
corresponding to the eigenvalue of $\mu = (1.6-3.0\bi)\tunit^{-1}$, the maximum LAP index is 137, 
corresponding to the length scale of $4.6\lunit$.
\textbf{Right:} 
$\vert \bg \vert = 0.3745\qunit$, $\delta = \Delta = 25 \tunit$, $b=104.2\bunit$, the maximum BT index is 752, 
corresponding to the eigenvalue of $\mu = (8.9+3.3\bi)\tunit^{-1}$, the maximum LAP index is 347, 
corresponding to the length scale of $2.4\lunit$.
The geometry is the neuron {\it 03b\_spindle4aACC}.  The gradient direction is 
$\bug =[0.7071  \ \ 0  \ \  0.7071]$.  
}
\end{figure}

To provide an at-a-glance summary of the significant BT and Laplace eigenfunctions at the end of 
the first gradient pulse $t=\delta$, we computed the matrix  
\ben
A(\delta) = \begin{bmatrix}
V_{11}e^{-\mu_1 \delta } &  &  & & \\
 &V_{12}e^{-\mu_2 \delta } &   &  &\\
 & &  & \ddots &  \\
&  & &  & V_{1N_{eig}}e^{-\mu_{N_{eig}} \delta }
\end{bmatrix}V^{-1}, \quad A(\delta) \in \C^{N_{eig}\times N_{eig}},
\een
for the neuron {\it 03b\_spindle4aACC} at two gradient amplitudes and two values of $\delta$.
It is easy to see that the magnetization solution at $t=\delta$ is
\ben
M(\bx,\delta) = \rho \sqrt{\vert \Omega \vert }
\left (V_{11}e^{-\mu_1 \delta }\psi_1(\bx)+\cdots+V_{1N_{eig}}e^{-\mu_{N_{eig}} \delta }\psi_{N_{eig}}(\bx)\right),
\een
where we time evolved the BT eigenfunctions according to their eigenvalues.
One can project $M(\bx,\delta)$ onto the space of the Laplace eigenfunctions:
\be{eq:lap_coeff_tsdelta}
M(\bx,\delta) = \rho \sqrt{\vert \Omega \vert } \left( c_1(\delta)\phi_1(\bx)+\cdots+c_{N_{eig}}(\delta)\phi_{N_{eig}}(\bx)\right).
\ee
The sum of $A(\delta)$ along the $k$th column is exactly $c_k(\delta)$. 
In Figure \ref{fig:lapbt_grid}, we show the indices $j$ and $k$ such that 
$\vert A_{jk}(\delta)\vert \geq 0.001$ for three choices of gradient amplitudes and pulse durations. 
This threshold was made so that the resulting signal is within 1\% of the signal where no entries 
of $A(\delta)$ were zeroed out.  
The row index of $A(\delta)$ is the index of the BT eigenfunction, the column index of $A(\delta)$ is 
the index of the Laplace eigenfunction.  The value of the $j$th row and the $k$th column, $A_{jk}(\delta)$,
is the projection onto the $k$th Laplace eigenfunction of the time evolved $j$th BT eigenfunction, $V_{1j}e^{-\mu_j \delta}\psi_j(\bx)$.  
We see a clear reduction of the number of significant BT eigenfunctions as $\delta$ increases.  At the 
higher gradient amplitude, there are a larger number of significant BT eigenfunctions as well as a larger number of Laplace eigenfunctions.  At $\delta=2.5\tunit$, the minimum significant Laplace eigenfunction length scale went from $3.5\lunit$ to $2.4\lunit$ as $\vert \bg \vert$ was increased from $0.075\qunit$ to $0.3745\qunit$.
Finally, we caution that $\vert c_k(\delta)\vert = \vert \sum_j  A_{jk}(\delta) \vert $ can be significantly smaller than $\sum_j \vert A_{jk}(\delta)\vert$ due to cancellations of numbers of opposite signs.  This means 
non-zero columns in Figure \ref{fig:lapbt_grid} may have a coefficient $c_k(\delta)$ that is negligible.

\subsection{Simulation of the apparent diffusion coefficient}
%\marginparnew{This Section is new}

As an additional validation of our computations of the Matrix Formalism signal, we compare the ADC
from the Matrix Formalism signal with several other computation methods.  The first method 
is a linear fit of the BTPDE signal at two small $b$-values, which we chose to be $b=0\bunit$ and $b=1\bunit$.
The second method uses the HADC model\cite{Haddar2018} which is derived from the BTPDE using 
mathematical homogenization.  The HADC model is a PDE model whose solution gives the time-dependent and gradient direction-dependent effective diffusion coefficient for arbitrary diffusion encoding sequences.
The third method is an analytical short time approximation formula for the ADC.
A well-known formula for the ADC in the short diffusion time regime is the following
short time approximation (STA) \cite{Mitra1992,Mitra1993}:
\ben
STA = \Dintr \left( 1 -  \frac{4\sqrt{\Dintr}}{3 \; \sqrt{\pi}}\sqrt{\Delta} \frac{{A}}{dim \;V} \right),
\een
where $\dfrac{A}{V}$ is the surface to volume ratio
and $\Dintr$ is the intrinsic diffusivity coefficient.  In the above formula the 
pulse duration $\delta$ is assumed to be very small compared to $\Delta$. 
A recent correction to the above formula \cite{Haddar2018}, taking into account the 
finite pulse duration $\delta$ and the gradient direction $\bug$, is the following:
\be{eq:STA_FP}
STA = \Dintr  \left[ 1- \frac{4 \sqrt{\Dintr}}{3 \; \sqrt{\pi} } 
					C_{\delta,\Delta}
					 \frac{A_\bug}{V } \right],
\ee
where 
\ben
A_\bug = \int_{\partial \Omega} \left(\bug \cdot \bn\right)^2\,ds
\een
is gradient direction dependent, and
\begin{align*} 
C_{\delta,\Delta} &= \dfrac{4}{35} 
					 \dfrac{\left(\Delta+\delta \right)^{7/2}+ \left(\Delta - \delta \right)^{7/2} - 2 \left( \delta^{7/2}+\Delta^{7/2} \right)  }{\delta^2 \left(\Delta-\delta/3\right)}  \label{eqn_def_C_delDEL}
				 = \sqrt{\Delta} \left( 1 + \dfrac{1}{3} \dfrac{\delta}{\Delta} - \dfrac{8}{35} \left( \dfrac{\delta}{\Delta}\right)^{3/2} + \cdots  \right).
\end{align*}
When $\delta \ll \Delta$, the value $C_{\delta, \Delta}$ is approximately $\sqrt{\Delta}$.

\begin{figure}[ht!]
 \centering
\includegraphics[width=0.5\textwidth]{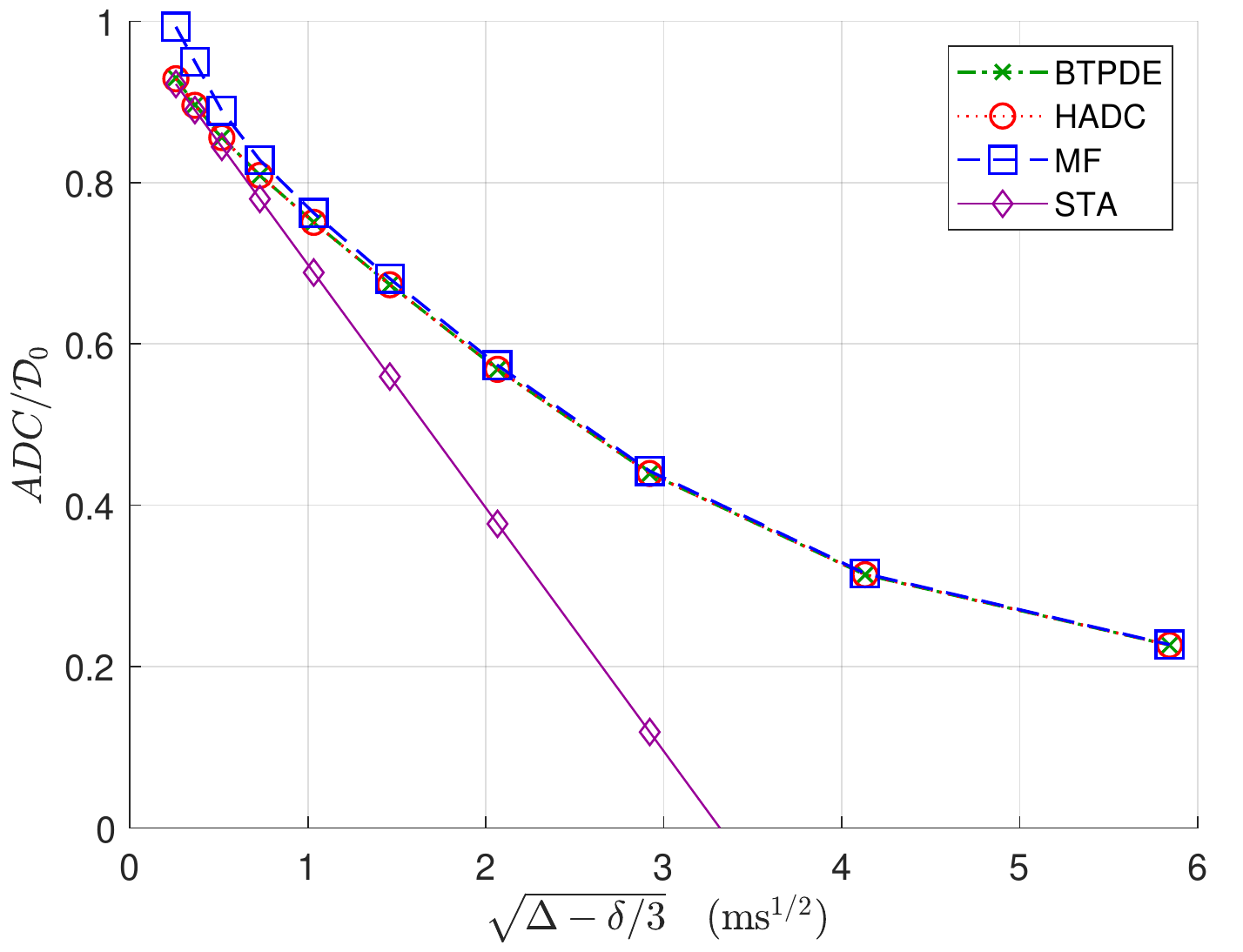}\\
  \caption{\label{fig:sta}  The computed apparent diffusion coefficient for PGSE sequences 
with $\Delta=2\delta$, where $\delta$ = $0.04\tunit$, $0.08\tunit$, $0.16\tunit$, $0.32\tunit$,
$0.64\tunit$, $1.28\tunit$, $2.56\tunit$, $5.12\tunit$, $10.24\tunit$, $20.48\tunit$.
The diffusion time is defined as $\Delta-\delta/3$.
The BTPDE ADC is obtained uses a linear fit of the ADC from $b=0$ and $b=1\bunit$.
The HADC model is a PDE model whose solution gives the time-dependent and gradient direction-dependent 
effective diffusion coefficient for arbitrary diffusion encoding sequences and it is derived from the BTPDE.
The STA is an analytical formula for the time-dependent and gradient direction-dependent effective
diffusion coefficient based on the surface area and the volume.
The MF is the Matrix Formalism ADC using $l_s^{min} = 0.5 \lunit$.  
The BTPDE, the HADC, and the MF simulations are computed on a 
finite element mesh of (44201 nodes, 160205 elements), with ODE solver tolerances $atol = 10^{-6}$,
$rtol=10^{-4}$.  For MF, this resulted in 9093 Laplace eigenfunctions.
The gradient direction is $\bug =[1  \ \ 0  \ \  0]$.  
The geometry is the neuron {\it 03b\_spindle4aACC}.  The intrinsic diffusion coefficient is 
$\Dintr = 2\times 10^{-3}\dunit$.
}
\end{figure}

We compute the Matrix Formalism ADC for the neuron {\it 03b\_spindle4aACC} 
using $l_s^{min} = 0.5 \lunit$, on a 
finite element mesh with 44201 nodes and 160205 elements, this resulted in 9093 Laplace eigenfunctions.
The BTPDE and the HADC simulations are also computed on this finite elements mesh, 
with ODE solver tolerances $atol = 10^{-6}$, $rtol=10^{-4}$.  
Ten PGSE sequences were simulated 
with $\Delta=2\delta$, where $\delta$ = $0.04\tunit$, $0.08\tunit$, $0.16\tunit$, $0.32\tunit$,
$0.64\tunit$, $1.28\tunit$, $2.56\tunit$, $5.12\tunit$, $10.24\tunit$, $20.48\tunit$.
In Figure \ref{fig:sta}, we see that the ADC from BTPDE and HADC are identical for all
the 10 PGSE sequences.  We see the STA matches them at the three shortest diffusion times, 
and the MF does not match them at the shortest 4 diffusion times, but starts to match them 
from the 5th diffusion time onwards.  These results are not surprising given that the required number of 
Laplace eigenfunctions to accurately represent the diffusion Green's function goes to infinity
as the diffusion time goes to 0.

\section{Discussion}

We have shown some of the functionalities of the Matrix Formalism Module 
within the diffusion MRI simulation toolbox SpinDoctor.  
We showed that the numerically computed $S^{\text{MF}}$ is very close to the reference signal from the 
Bloch-Torrey PDE for realistic neuron geometries at a wide range of b-values, 
and the agreement between the two are good at a wide range of diffusion times.  
We examined in detail the contributions of different eigenmodes at different diffusion times and different b-values.
By choosing to represent the eigenvalues by a quantity of length $l_s(\lambda_n)$, we highlighted the important spatial scales 
that contribute to the diffusion MRI signal. 
\soutnew{}{It can be seen that 
the longest length scales of the Laplace eigenfunctions roughly relate to the size of the geometry
in the spatial directions.
In addition, the length scale of a Laplace eigenfunction is related to the "wavelength" of the oscillations of the eigenfunction.}

\soutnew{}{We have validated our Matrix Formalism computations by comparing it with the reference 
solution obtained from solving the BTPDE for a whole neuron.  
We showed that the $S^{\text{MF}}$ can be made more 
accurate by decreasing the requested minimum length scale $l_s^{min}$ and refining the finite elements mesh. 
We also showed through
an ADC computation that, fixing the number of Laplace eigenfunctions a priori means that $S^{\text{MF}}$ 
cannot be made accurate at arbitrarily short diffusion times, since the number of required Laplace eigenfunctions, to describe
the magnetization limit as the diffusion time goes to 0, goes to infinity.}

\soutnew{}{We gave the transformation that 
link the Laplace eigenfunctions to the eigenfunctions of the Bloch-Torrey operator.  We computed the BT eigenfunctions and eigenvalues for a whole neuron.  We showed that at high enough gradient amplitude, 
the first few BT eigenfunctions (ordered by their real part) have support that are localized to 
a small region of the geometry.  We showed the BT eigenvalues have a larger range in both their real
part and imaginary part as the gradient amplitude increases.  This explains the faster transient behavior and 
the oscillatory nature of the magnetization at high gradient amplitude.  We showed that the higher indexed
BT eigenfunctions have contributions from Laplace eigenfunctions that have shorter length scales.}

There are two important advantages to the Matrix Formalism signal representation \soutnew{}{over the other
representations}.  The first advantage is analytical, 
this representation makes explicit 
the link between the Laplace eigenvalues and eigenfunctions of the biological cell and its diffusion MRI signal.  
This clear link may help in the formulation of reduced models of the diffusion MRI signal that is closer to the physics of the problem.  The second advantage is numerical, once the Laplace eigendecomposition has been computed and saved, 
the diffusion MRI signal can be calculated for arbitrary diffusion-encoding sequences and b-values 
at negligible additional cost.  This will make it possible to use the Matrix Formalism as the inner loop of optimization procedures. 

The need for a mathematically rigorous model of the diffusion MRI signal arising from realistic cellular structures was re-iterated in recent review papers \cite{Novikov2019,Novikov2018,Fieremans2018}. 
Given that \soutnew{Bloch-Torrey equation}{the Bloch-Torrey equation posed in realistic geometries} is a gold-standard reference model of the diffusion MRI signal and 
the Matrix Formalism signal representation is equivalent to the reference model as long as enough eigenmodes 
are included, Matrix Formalism may be a possible bridge to formulating practical "inverse models" that can
be used to robustly estimate biological relevant parameters from the acquired experimental data.  
In this paper, we moved Matrix Formalism a step closer to being a practically computable model and 
showed the number of significant eigenmodes is around \soutnew{100}{a few hundred} for realistic neurons.  
The next step may be searching for a unique set of "modes" onto which to project the eigenmodes 
of a population of many neurons in a voxel.  Finding such a universal set of "modes" would require
more advanced mathematical analysis on the diffusion operator in geometries with multiple length 
scales.  A modified Fourier basis may be considered such as in \cite{greengard1990} as a possible future direction of research.

Currently, the Matrix Formalism Module allows the computation of the Matrix Formalism signal and the 
Matrix Formalism Gaussian Approximation signal for realistic neuron (impermeable membranes) 
with the PGSE sequence.
Matrix Formalism for permeable membranes 
and for general diffusion-encoding sequences are under development and will be released in the future.    
The SpinDoctor toolbox and the Neuron Module have been developed in the MATLAB R2017b and require no additional MATLAB toolboxes.  However, the current version of the Matrix Formalism Module requires the
MATLAB PDE Toolbox (2017 or later) due to certain difficulties of implementing the matrix eigenvalue solution
on a restricted eigenvalue interval.  This technical issue will be addressed in a future release.      

The Matrix Formalism Module follows the same workflow as SpinDoctor and the Neuron Module 
and builds upon the functionalities of SpinDoctor. 
To use the Matrix Formalism Module, it is necessary to read first the documentation of SpinDoctor.
The source code, examples, documentation of SpinDoctor, the Neuron Module and the Matrix Formalism Module 
are available at \url{https://github.com/jingrebeccali/SpinDoctor}.

\soutnew{In the Appendix, we list the input files, as well as important quantities and functions relevant to the Matrix Formalism Module, noting where relevant, the input parameters that are not applicable (marked by "na")
to the current version of the Matrix Formalism Module.  Sample output figures are also provided.
}{}
\section{Conclusion}
We presented a simulation module that we have implemented inside a MATLAB-based diffusion MRI simulator called SpinDoctor that efficiently computes the Matrix Formalism representation for realistic geometrical models of neurons.
With this new simulation tool, we seek to bridge the gap between physical quantities
closely related to the cellular geometrical structure, namely, its Laplace eigenfunctions, eigenvalues and their length scales,
with the measured diffusion MRI signal. 
We hope this Matrix Formalism Module makes the mathematically rigorous signal representation 
into a practical model for the research community.

\section*{Acknowledgment}
The authors gratefully acknowledge the {\it French-Vietnamese Master in Applied Mathematics} program 
whose students (co-authors on this paper,  {Try Nguyen Tran} and {Van-Dang Nguyen}) 
have contributed to the SpinDoctor project during their internships in France in the past 
several years, as well as the {\it Vice-Presidency for Marketing and International Relations} 
at Ecole Polytechnique for financially supporting a part of the students' stay.

\appendix
\section{Derivation of the Matrix Formalism signal representation}
\label{sec:appendix1}

Let $\phi_n(\bx)$ and $\lambda_n$, $n = 1,\cdots$,  be the $L^2$-normalized eigenfunctions and eigenvalues associated to the Laplace operator 
with homogeneous Neumann boundary conditions (the surface of the neurons are supposed
impermeable):
\begin{alignat*}{3}
-\nabla \Dintr \left(\nabla \phi_n(\bx)\right) &= \lambda_n \phi_n(\bx),&\quad \bx \in \Omega,\\
\Dintr \nabla \phi_n(\bx) \cdot \bnu(\bx) &= 0, &\quad \bx \in \Gamma.
\end{alignat*}
We assume the non-negative real-valued eigenvalues are ordered in non-decreasing order:
\ben
0 = \lambda_1 < \lambda_2 \leq \lambda_3 \leq \cdots
\een
so that $\lambda_1 = 0$ (this means the first Laplace eigenfunction is the constant function).
There is only one constant eigenfunction because we assume the neuron is a connected domain.
Let $L$ be the diagonal matrix containing the first $N_{eig}$ Laplace eigenvalues:
\be{}
L = \text{diag} [\lambda_1,\lambda_2,\cdots, \lambda_{N_{eig}}] \in \R^{N_{eig}\times N_{eig}}.
\ee

For the PGSE sequence, the time interval $[0,T_E], T_E =\delta+\Delta $ is separated into three subintervals:
\be{eq:pgse}
f(t) =
\begin{cases}
1, \quad & t \in [0, \delta], \\
0,\quad & t \in [\delta, \Delta],\\
-1,\quad & t \in [\Delta, T_E].
\end{cases}
\ee

If $t \in [0, \delta]$, decompose the solution $M(\bx, t)$ in the basis of $\{ \phi_n (\bx) \}$ as
\begin{equation}
M(\bx, t) = \sum _{n = 1} ^{N_{eig}} T_n(t) \phi_n (\bx) =  [\phi (\bx) ]^T T(t),
\end{equation}
where we define the space and time dependent column vectors 
\ben
\phi (\bx) = \begin{bmatrix}  \phi_1(\bx)  \\ \vdots  \\ \phi_{N_{eig}}(\bx)   \end{bmatrix}, \ 
T(t) = \begin{bmatrix} T_1 (t)  \\ \vdots \\ T_{N_{eig}}(t) \end{bmatrix}. 
\een

Substituting into the original Bloch-Torrey equation, multiplying both sides with $\{ \phi_k (\bx) \}$ and integrating over $\Omega$ give
\begin{equation} \label{eq:BTsystem}
\frac{\partial}{\partial t}{T_k(t)} = -\lambda_k T_k(t) - \bi\gamma \sum_ {n = 1}^{N_{eig}} T_n(t)\int_{\Omega} \bg \cdot \bx \phi_n(\bx)\phi_k(\bx) d\bx, \quad k = 1,..., N_{eig}.
\end{equation}

Let $W(\bg)$ be the ${N_{eig}\times N_{eig}}$ matrix whose entries are:
\be{}
W(\bg) =  \left[ \int_{\Omega} \bg \cdot \bx \phi_m(\bx)\phi_n(\bx) d\bx \right] =g_x A^x+g_y A^y + g_z A^z
\ee
where $\bg = (g_x,g_y,g_z),$ and
\be{}
A^{i} =\begin{bmatrix}
a_{mn}^{i}
\end{bmatrix}, \quad i = \{x,y,z\}, \quad 1 \leq m,n \leq N_{eig},
\ee
are three symmetric ${N_{eig}\times N_{eig}}$ matrices whose entries are 
the first order moments in the coordinate directions of the product of pairs of eigenfunctions:
\ben
\begin{split}
a_{mn}^x &=  \int_{\Omega}  x\phi_m(\bx) \phi_n(\bx) d\bx,\;\\
a_{mn}^y &=  \int_{\Omega}  y\phi_m(\bx) \phi_n(\bx) d\bx,\\
a_{mn}^z &=  \int_{\Omega}  z\phi_m(\bx) \phi_n(\bx) d\bx.
\end{split}
\een

We define the complex-valued matrix $K(\bg)$ 
\be{}
K(\bg) \equiv L + \bi \gamma W(\bg) .
\ee

Then the system of equations \ref{eq:BTsystem} can be written as 
\begin{equation}
\frac{\partial}{\partial t}{T(t)} = -K T(t).
\end{equation}
The solution to this system of linear differential equations is
\begin{equation}
T(t) = e^{-Kt} T(0), \ t \in [0, \delta].
\end{equation}

The initial condition
\ben
M(\bx,0) = \sum _{n = 1} ^{N_{eig}} T_n(0) \phi_n (\bx) = [\phi (\bx) ]^T T(0),
\een
gives us the initial coefficients
\begin{equation}
\ T_n (0) = \int_\Omega M(\bx,0) \phi_n(\bx) d \bx.
\end{equation}
We get
\begin{equation}
M(\bx, t) = [\phi (\bx) ]^T e^{-K t} T(0), \ t \in [0,\delta].
\end{equation}

Similar analysis can be applied to the interval $[\delta,  \Delta]$ with $f(t) = 0$ and $M(\bx,\delta)$ acting as the initial condition  for the system of differential equations. At the end of this interval, the magnetization is represented as
\begin{equation}
M(\bx, \Delta) = [\phi (\bx) ]^T e^{-L (\Delta-\delta)} e^{-K \delta} T(0).
\end{equation}

The magnetization measured at the end of the second pulse is
\begin{equation} \label{M_TE}
M(\bx, T_E) = [\phi (\bx) ]^T e^{-K^* \delta} e^{-L (\Delta-\delta)} e^{-K \delta} T(0).
\end{equation}

Denote $H(\bg, f) = e^{-K^* \delta} e^{-L (\Delta-\delta)} e^{-K \delta} $ and $\Phi = \begin{bmatrix} \int_\Omega  \phi_1(\bx) d\bx \\ \vdots  \\ \int_\Omega \phi_{N_{eig}}(\bx) d\bx  \end{bmatrix}$. Then the echo signal is computed by integrating $M(\bx, T_E)$ over $\Omega$: 
\begin{equation} \label{MF_signal}
S = \int_\Omega M(\bx, T_E) d \bx = \Phi^T H(\bg, f) T(0).
\end{equation}

If $M(\bx,0) = \rho$, where $\rho$ is the initial spin density, then $T(0) = \rho \Phi$. With homogeneous Neumann boundary conditions, the $L^2$-normalized eigenfunctions $\phi_n(\bx)$ have the following property 
\be{eq:phi}
 \Phi_n = \int_ \Omega \phi_n(\bx)  d\bx =
\begin{cases}
\sqrt{|\Omega|}, \quad & n =1, \\
0,\quad & n \geq 2,
\end{cases}
\ee
with, in particular, $\phi_1(\bx) \equiv \sqrt{\frac{1}{\vert \Omega \vert }}$.
The vector $\Phi$ has only one non-zero value in the first position. Eq. \ref{MF_signal} is then reduced to 
\begin{equation} \label{sig_TE}
S = \rho |\Omega| H_{11} (\bg, f) = S_0 H_{11} (\bg, f).
\end{equation}
Note that since $K$ and $K^*$ are symmetric matrices, $e^{-K^* \delta}$ and $e^{-K \delta}$ are also symmetric; thus, we have:
\begin{equation}
H^T(\bg, f) =   \left[e^{-K \delta} \right]^T e^{-L (\Delta-\delta)} \left[e^{-K^* \delta}\right]^T = e^{-K \delta}  e^{-L (\Delta-\delta)} e^{-K^* \delta}
\end{equation}
Eq. \ref{M_TE} can be written as 
\begin{equation}
M(\bx, T_E) = [ T(0) ]^T H^T(\bg, f) \phi (\bx)=[ T(0) ]^T e^{-K \delta} e^{-L (\Delta-\delta)} e^{-K^* \delta} \phi (\bx).
\end{equation}
In Eq. \ref{sig_TE} we can reverse the order of the matrix exponentials in $H(\bg, f)$ and achieve the same echo signal:
\begin{equation}
S =  S_0 \left[ e^{-K \delta} e^{-L (\Delta-\delta)} e^{-K^* \delta} \right] _{11}.
\end{equation}

\section{Eigenvalues and eigenfunctions of the Bloch-Torrey operator}
\label{sec:appendix2}
To relate the eigen-decomposition of the Laplace operator $-\Dintr \nabla ^2$ to the eigen-decomposition of the 
complex-valued Bloch-Torrey operator $-\left(\Dintr \nabla ^2 - \bi\gamma  \bg \cdot \bx\right)$, we first 
define $\psi_m(\bx)$ and $\mu_m$, $m = 1,\cdots, N_{eig}$, to be the eigenfunctions and eigenvalues associated to the Bloch-Torrey operator subject to homogeneous Neumann boundary conditions:
\begin{alignat}{3}
\label{eq:BT_operator}
-(\Dintr \nabla ^2 - \bi \gamma  \bg \cdot \bx) \psi_m(\bx) &= \mu_m \psi_m(\bx), &\quad \bx \in \Omega,   \\
\Dintr \nabla \psi_m(\bx) \cdot \bnu(\bx) &= 0, &\quad \bx \in \Gamma. 
\end{alignat}
Representing $\psi_m(\bx), m = 1,\cdots, N_{eig}$, in the basis of the Laplace eigenfunctions
$\phi_n(\bx), n = 1,\cdots, N_{eig}$:
\be{}
\psi_m (\bx) = \sum_{n = 1} ^{N_{eig}} c_{mn} \phi_n (\bx) , \ m = 1,...,N_{eig},
\ee 
and putting it into Eq. \ref{eq:BT_operator}, 
\begin{equation}
- \sum_ {n = 1}^{N_{eig}}c_{mn}\Dintr \nabla ^2 \phi_n(\bx) + \bi\gamma  \bg \cdot \bx \sum_ {n = 1}^{N_{eig}}c_{mn}\phi_n(\bx) = \mu_m \sum_ {n = 1}^{N_{eig}}c_{mn}\phi_n(\bx),
\end{equation}
we get
\begin{equation}
 \sum_ {n = 1}^{N_{eig}}c_{mn} \lambda_n \phi_n(\bx) + \bi\gamma  \bg \cdot \bx \sum_ {n = 1}^{N_{eig}}c_{mn}\phi_n(\bx) = \mu_m \sum_ {n = 1}^{N_{eig}}c_{mn}\phi_n(\bx).
\end{equation}
Multiply both sides with $\phi_k(\bx)$ and integrate over $\Omega$:
\begin{equation}
\label{eq:mateq}
\lambda_k c_{mk}+ \bi\gamma \sum_ {n = 1}^{N_{eig}}c_{mn} \int_{\Omega} \bg \cdot \bx \phi_n(\bx)\phi_k(\bx)d\bx = \mu_m c_{mk}, \quad m,k  = 1,..., N_{eig}.
\end{equation}
Then writing Eq. \ref{eq:mateq} as a matrix equation for $C=[c_{mk}]\in \C^{N_{eig}\times N_{eig}}$, where $m$ is the row index and $k$ is the column index:
\ben
C \left(L + \bi\gamma W(\bg)\right) = CK(\bg) = \Sigma C.
\een
We see that we can obtain $C$ and $\Sigma$ by diagonalizing the complex-valued matrix $K(\bg)$:
\be{}
K(\bg) = V \Sigma V^{-1},
\ee
where the matrix $V$ has the eigenvectors in the columns
and $\Sigma$ has the eigenvalues on the diagonal.
We can obtain $C$ by setting it to  
\be{}
C \equiv V^{-1}.
\ee 
The diagonalization of the complex conjugate of $K(\bg)$ is  
\be{}
\left(K(\bg)\right)^{*} = L - \bi \gamma W(\bg) = (V^{-1})^{*} \Sigma^{*} (V)^{*}.
\ee
At the end of the first gradient pulse, the magnetization is
\begin{equation}
M(\bx,\delta) = \mathcal M_0^{\psi}e^{-\Sigma \delta } [\psi_1(\bx),\cdots\psi_{N_{eig}}(\bx)]^T,
\end{equation}
where $\mathcal M _0^{\psi}\in \C^{1 \times N_{eig}}$ is the vector of coefficients of the 
initial magnetization projected onto the eigenfunctions of $K(\bg)$.
Similarly, at the beginning of the second time interval, the magnetization is
\begin{equation}
M(\bx,\Delta) = \mathcal M_\delta^{\phi} e^{-L (\Delta-\delta) }[ \phi_1(\bx),\cdots\phi_{N_{eig}}(\bx)]^T,
\end{equation}
where $\mathcal M_\delta^{\phi}\in \C^{1\times N_{eig}}$ is the vector of coefficients of the 
magnetization at $t=\delta$ projected onto the eigenfunctions of $L$.
Finally, at the end of the second interval, the magnetization is 
\begin{equation}
M(\bx,\delta+\Delta) = \mathcal M_\Delta^{\tilde{\psi}} e^{- \Sigma^{*} \delta } [\tilde{\psi}_1(\bx),\cdots \tilde{\psi}_{N_{eig}}(\bx)]^T,
\end{equation}
where $\mathcal M_\Delta^{{\tilde{\psi}}}\in \C^{1\times N_{eig}}$ is the vector of coefficients of the 
magnetization at $t=\Delta$ projected onto the eigenfunctions of $K(\bg)^{*}$ (we denoted this $\tilde{\psi}$).
The conversion between $\psi$ and $\phi$ is 
\be{}
\begin{split}
[ \psi_1(\bx),\cdots\psi_{N_{eig}}(\bx)]^T &= V^{-1} [ \phi_1(\bx),\cdots\phi_{N_{eig}}(\bx)]^T \\
V [ \psi_1(\bx),\cdots\psi_{N_{eig}}(\bx)]^T &= [ \phi_1(\bx),\cdots\phi_{N_{eig}}(\bx)]^T 
\end{split}.
\ee
So 
\be{}
\mathcal M_0^{\psi} = \rho \sqrt{\vert \Omega \vert }[V_{11},\cdots,V_{1n}],
\ee
is the first row of $V$, and 
\be{}
\mathcal M_\delta^{\phi} = \mathcal M_0^{\psi}e^{-\Sigma\delta }  V^{-1}.
\ee
The conversion between $\tilde{\psi}$ and $\phi$ is 
\be{}
\begin{split}
[ \tilde{\psi}_1(\bx),\cdots \tilde{\psi}_{N_{eig}}(\bx)]^T &= V^{*} [ \phi_1(\bx),\cdots\phi_{N_{eig}}(\bx)]^T \\
(V^{-1})^{*} [ \tilde{\psi}_1(\bx),\cdots \tilde{\psi}_{N_{eig}}(\bx)]^T &= [ \phi_1(\bx),\cdots\phi_{N_{eig}}(\bx)]^T 
\end{split}
\ee
Thus, 
\be{}
\mathcal M_\Delta^{{\tilde{\psi}}} =\mathcal M_\delta^{\phi} e^{-L (\Delta-\delta) } (V^{-1})^{*}.
\ee
The final magnetization is therefore
\be{}
M(\bx,\delta+\Delta) =  \rho \sqrt{\vert \Omega \vert }[V_{11},\cdots,V_{1n}] e^{-\Sigma \delta }  V^{-1}  e^{-L (\Delta-\delta) } (V^{-1})^{*} e^{- (\Sigma)^{*} \delta }V^{*} [ \phi_1(\bx),\cdots\phi_{N_{eig}}(\bx)]^T.
\ee
Since the first Laplacian eigenfunction is the constant function and it is orthogonal to all the other 
Laplacian eigenfunctions, the signal is 
\be{}
S = \int_{\bx\in \Omega} M(\bx,\delta+\Delta)d\bx =  \rho {\vert \Omega \vert }[V_{11},\cdots,V_{1n}] e^{-\Sigma\delta }  V^{-1}  e^{-L (\Delta-\delta) } (V^{-1})^{*} e^{- (\Sigma)^{*} \delta }[V_{11},\cdots,V_{1n}]^{*},
\ee
where we only kept the first column of $V^{*}$.

\bibliography{myref_final}%

\end{document}